\documentclass[preprint,12pt]{elsarticle}
\usepackage{amsmath}
\usepackage{amssymb}
\usepackage{graphicx}
\usepackage{hyperref}
\usepackage{color}
\usepackage{multicol}
\usepackage{calrsfs}
\usepackage{algorithm}
\usepackage{algorithmic}
\usepackage{cancel}
\usepackage{multirow}
\usepackage[algo2e]{algorithm2e} 
\usepackage{graphicx,epsfig}
\usepackage{amsmath,amssymb,euscript}
\usepackage{latexsym}
\usepackage{color}
\newtheorem{theorem}{\bf Theorem}[section]

\newtheorem{definition}[theorem]{\bf Definition}

\newtheorem{proposition}[theorem]{\bf Proposition}

\newtheorem{corollary}[theorem]{\bf Corollary}
\newtheorem{remark}[theorem]{\bf Remark}

\newcommand{\Proofend}{\hfill$\diamondsuit$}
\newcommand{\PP}{{\mathbb P}}
\newcommand{\RR}{{\mathbb R}}

\newcommand{\NN}{\mathbb{N}}

\newcommand{\C}{\mathcal{C}}

\newcommand{\T}{\mathcal{T}}
\newcommand{\V}{\mathcal{V}}
\newcommand{\W}{\mathcal{W}}
\newcommand{\bigO}{\mathcal{O}}

\begin{document}
\begin{frontmatter}

\title{
Orthonormal polynomial wavelets associated with de la Vall\'ee Poussin-type interpolation on $[-1,1]$} 


\author[inst1]{Woula Themistoclakis
\corref{cor1}}
\ead{woula.themistoclakis@cnr.it}

\author[inst2]{Marc Van Barel}
\ead{marc.vanbarel@cs.kuleuven.be}

\cortext[cor1]{Corresponding author.}

\affiliation[inst1]{organization={CNR National
        Research Council of Italy,
        IAC Institute for Applied Computing ``Mauro Picone''},
            addressline={via P. Castellino, 111}, 
            city={Napoli}, 
            postcode={80131}, 
            country={Italy}}

\affiliation[inst2]{organization={Department of Computer Science, KU Leuven},
            addressline={Celestijnenlaan 200A}, 
            city={Leuven (Heverlee)}, 
            postcode={B-3001}, 
            country={Belgium}}

\begin{abstract}
Starting from de la Vall\'ee Poussin type (VP) interpolation, the authors have recently introduced a family of interpolating polynomial scaling and wavelet bases generating the approximation and detail spaces of a non-standard multiresolution analysis. Motivated by the fact that, in many applications, orthonormal rather than interpolating bases are preferable, the present study develops a new family of scaling and wavelet polynomials that provide well-localized and orthonormal bases for the same approximation and detail spaces.

We show that the proposed new bases have a behavior very similar to the interpolating bases already introduced, presenting similar features although they are not interpolating but orthonormal. In particular, we study the Fourier projection corresponding to the proposed orthonormal scaling basis, 
and introduce a discrete version of it by approximating the Fourier--like coefficients. For both continuous and discrete orthogonal projections, we prove the uniform boundedness of the 
Lebesgue constants and the uniform convergence with an asymptotic rate comparable  with the best uniform polynomial approximation. 

Numerical experiments confirm the theoretical results and compare the new orthonormal VP scaling and wavelet bases with the interpolating case previously treated by the authors.
\end{abstract}


\end{frontmatter}

\section{Introduction}
In \cite{TB-wave} the authors recently introduced a family of interpolating polynomial wavelets $\psi_{n,k}^m(x)$, $k=1,\ldots, 2n$, which, at each resolution level $n\in\NN$, depend on a free parameter $m=\lfloor \theta n\rfloor$, with $\theta\in ]0,1[$ arbitrarily fixed. 

Like other polynomial wavelets (see, e.g., \cite{PSTasche, FisPre, FT-wave}) they are not generated by dilation and translation of a single mother function.  Their construction generalizes the trigonometric case considered in \cite{Prestin-trig}, and is based on de la Vall\'ee Poussin (VP, for short) interpolation at the Chebyshev nodes:
\[
X_n=\left\{x_k^n=\cos\frac{(2k-1)\pi}{2n},\ k=1,\ldots,n\right\}, \qquad n\in\NN .
\]
We recall that VP interpolation, originally introduced in \cite{Th-1999}, provides near-best polynomials for the uniform approximation of arbitrary continuous functions, reducing the Gibbs phenomenon and yielding optimal error estimates, even with respect to the $L^1$ norm (see \cite{Th-2012, Th-L1, OT-APNUM21, OT-DRNA21} and the references therein for more details).

The basis polynomials of this interpolation process are 
the fundamental VP polynomials $\Phi_{n,k}^m$,  $k=1,\ldots,n$, defined in \eqref{sca}. In \cite{TB-wave}, from the wavelet perspective, they are interpreted as  scaling functions at level $n$ with parameter $m$. They satisfy the interpolation property
\begin{equation}\label{sca-inter}
\Phi_{n,k}^m(x_h^n)=\delta_{k,h}, \qquad h,k=1,\ldots,n,
\end{equation}
and form a well--localized, Riesz--stable basis of the sampling space 
\begin{equation}\label{Sca-space}
\V_n^m:=\mathrm{span}\{\Phi_{n,k}^m \: \ k=1,\ldots,n\}.
\end{equation}
This is a polynomial space of dimension $n$ which is nested   between the canonical spaces  $\PP_\nu$ of polynomials of degree at most $\nu$, as follows
\begin{equation}\label{nested}
 \PP_{n-m}\subseteq \V_n^m\subset \PP_{n+m-1},\qquad 
 0<m<n. 
\end{equation}
As in classical multiresolution analysis, the resolution level $n$ can be increased or decreased; however, instead of the usual scaling factor 2, multiples of 3 are considered in \cite{TB-wave}.

More precisely, since $\V_n^m\subset\V_{3n}^m$, the detail space
$\W_n^m$ is defined as the orthogonal complement of $\V_n^m$ in $\V_{3n}^m$, i.e.,
\begin{equation}\label{compl}
\V^m_ {3n}=\V^m_{n}\oplus \W_n^m \quad\mbox{and}\quad \V_n^m\perp\W_n^m,
\end{equation}
where orthogonality is understood with respect to the following scalar product
\[
<f,g>_{L^2_w}:=\int_{-1}^1f(x)g(x)w(x)dx,\qquad w(x):=\frac 1{\sqrt{1-x^2}}.
\]
Moreover, since $X_n\subset X_{3n}$, by considering the set of nodes 
\[
Y_n:=X_{3n}-X_n=\{y_h^n\ :\ h=1,\ldots, 2n\},
\] 
a well--localized, Riesz--stable, and interpolatory wavelet basis of the detail space
\begin{equation}\label{Wav-space}
\W_n^m=\mathrm{span}\{\psi_{n,k}^m \: \ k=1,\ldots,2n\}
\end{equation}
is uniquely determined in \cite{TB-wave} by imposing the interpolation property
\begin{equation}\label{wave-inter}
 \psi_{n,k}^m(y_h^n)=\delta_{k,h}, \qquad h,k=1,\ldots,2n. \end{equation}
 Despite the non-orthogonality of the interpolating bases in \eqref{Sca-space} and \eqref{Wav-space}, efficient decomposition–reconstruction algorithms based on fast discrete cosine transforms were derived in \cite{TB-wave} and further generalized in \cite{Cotronei2025}.
In the latter, following a non-standard modeling \cite{ORT-Jmiv}, they were applied to image compression and denoising, achieving competitive performance compared to classical wavelets.
However, as explained in \cite{Cotronei2025}, the non-orthogonality of the interpolating scaling and wavelet bases necessitates the introduction of suitable normalization factors to control the loss of signal energy. In \cite{Cotronei2025}, these correction factors were theoretically derived for image compression and experimentally refined for signal denoising.
Nevertheless, the problem of determining optimal normalization factors, independent of the specific data and resolution levels, remains open.

In the present paper, we propose a different approach to this problem. For both the approximation and detail spaces at resolution level $n$, we define new, well-localized, orthonormal VP scaling and wavelet bases, dependent on the free parameter $m$. These bases behave similarly to the previous interpolating VP scaling and wavelet bases but, in addition, are orthonormal, with all the consequent advantages that this property entails.

The paper is organized as follows. Section~\ref{sec-scaling-functions}  introduces the orthonormal VP scaling functions. Section~\ref{sec-approx-oper} studies the Fourier projection and its discrete counterpart based on these scaling functions. Section~\ref{sec-wavelets} is devoted to the orthonormal VP wavelet functions. In Section~\ref{sec-dec-and-rec}, fast decomposition and reconstruction algorithms are developed. Section\ref{sec-experiments} presents numerical experiments comparing the continuous and discrete Fourier projections with the interpolating VP projection. Section~\ref{sec-proofs} contains the proofs of the theoretical results. Finally, Section~\ref{sec-conclusion} presents the conclusions of the present study and outlines directions for future research.
\section{Orthonormal VP scaling functions}\label{sec-scaling-functions}
\begin{figure}[!htb]%
\begin{center}
\includegraphics[scale = 0.40]{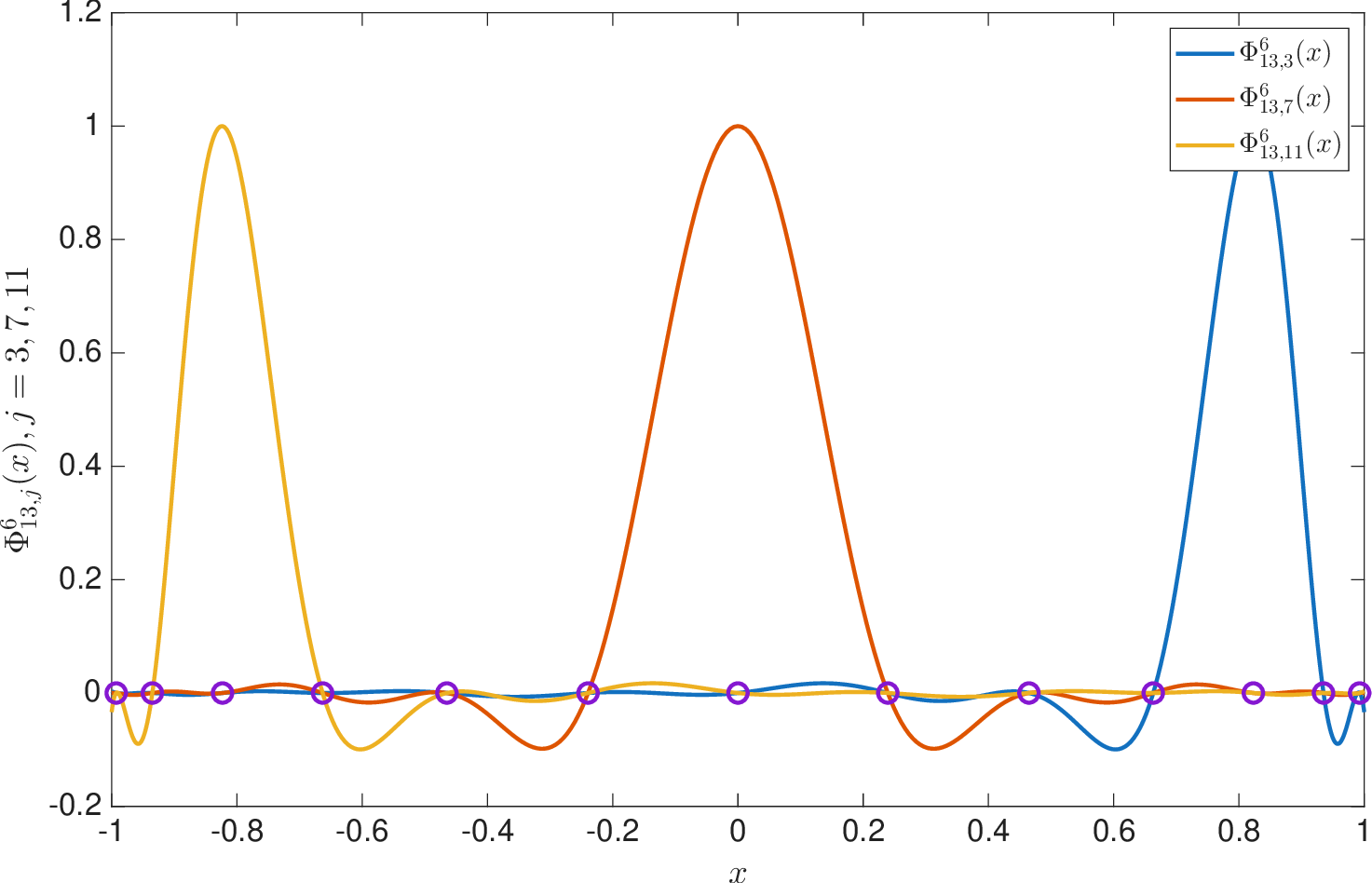}
\vspace{-.5cm}
\end{center}
\caption{ The interpolating VP scaling functions $\Phi_{n,k}^m$ for $n=13$, $m=6$, and $k=3, 7, 11$. The open dots denote the set of nodes $X_n$.\label{fig_Phi} }
\end{figure}

Let us first recall that on the set $X_n$ of the Chebyshev zeros 
the interpolating VP scaling functions of level $n$ and parameter $m<n$ are defined as follows \cite{TB-wave}
\begin{equation}\label{sca}
\Phi_{n,k}^m(x):=\frac\pi n\sum_{r=0}^{n+m-1}\mu_{n,r}^m p_r(x_k^n)p_r(x),\qquad k=1,\ldots,n\qquad x\in[-1,1],
\end{equation}
where 
\begin{equation}\label{muj}
\mu_{n,r}^m :=\left\{\begin{array}{ll}
1 & \mbox{if}\quad 0\le r\le n-m\\ [.1in]
\displaystyle\frac{m+n-r}{2m} & \mbox{if}\quad
n-m< r< n+m\\
0 & \mbox{otherwise}
\end{array}\right.
\end{equation}
and $p_n$ denotes the orthonormal Chebyshev polynomial of the first kind and degree $n$, namely
\begin{equation}\label{cheb-pol}
p_n(x):=\cos(n \arccos x)\left\{\begin{array}{ll}
\sqrt{ 1/\pi} &  \mbox{ if $n=0$}\\ [.1in]
\sqrt{2/\pi} &   \mbox{ if $n>0$}
\end{array}\right. \qquad |x|\le 1, \qquad n=0,1,\ldots
\end{equation}
In Figure~\ref{fig_Phi}, the polynomials $\Phi_{13,3}^6$, $\Phi_{13,7}^6$ and $\Phi_{13,11}^6$ are plotted. Note the interpolation property \eqref{sca-inter} at the Chebyshev zeros of degree $13$ marked with circles.


It is known that the scaling functions in \eqref{sca} constitute an interpolating Riesz--basis of the approximation space
\eqref{Sca-space}, which is also generated 
by the following rearrangement of the Chebyshev polynomials  
\begin{equation}\label{q-basis}
q_{n,r}^m(x):=\left\{\begin{array}{ll}
p_r(x)&\mbox{if}\quad 0\le r\le
n-m,\\ [.1in]
\displaystyle\mu_{n,r}^m p_r(x)-\mu_{n,2n-r}^m p_{2n-r}(x)
&\mbox{if}\quad n-m<r<n,
\end{array}\right.
\end{equation}
These polynomials constitute an orthogonal basis of $\V_n^m$ satisfying 
\begin{equation}\label{q-prod}
<q_{n,r}^m,\ q_{n,s}^m>_{L^2_w}=\delta_{r,s} \cdot \nu_{n,r}^m,  \qquad r,s=0,1,\ldots, n-1,
\end{equation}
\begin{equation}\label{nur}
\nu_{n,r}^m:=\left\{\begin{array}{ll} 1 & \mbox{if}\quad 0\le r\le
n-m,\\ [.1in] 
\displaystyle\frac{m^2+(n-r)^2}{2m^2} &
\mbox{if}\quad
 n-m<r< n,
\end{array}\right.\qquad r=0,1,\ldots,n-1.
\end{equation}
Moreover, we recall the following formula for the change of bases
\begin{equation}\label{change-sca}
\Phi_{n,k}^m(x)=\frac \pi n \sum_{r=0}^{n-1}p_r(x_k^n) q_{n,r}^m(x), \qquad k=1,2,\ldots,n,
\end{equation}
and the inverse formula given in \eqref{q-fi}.

\begin{figure}[!htb]%
\begin{center}
\includegraphics[scale = 0.40]{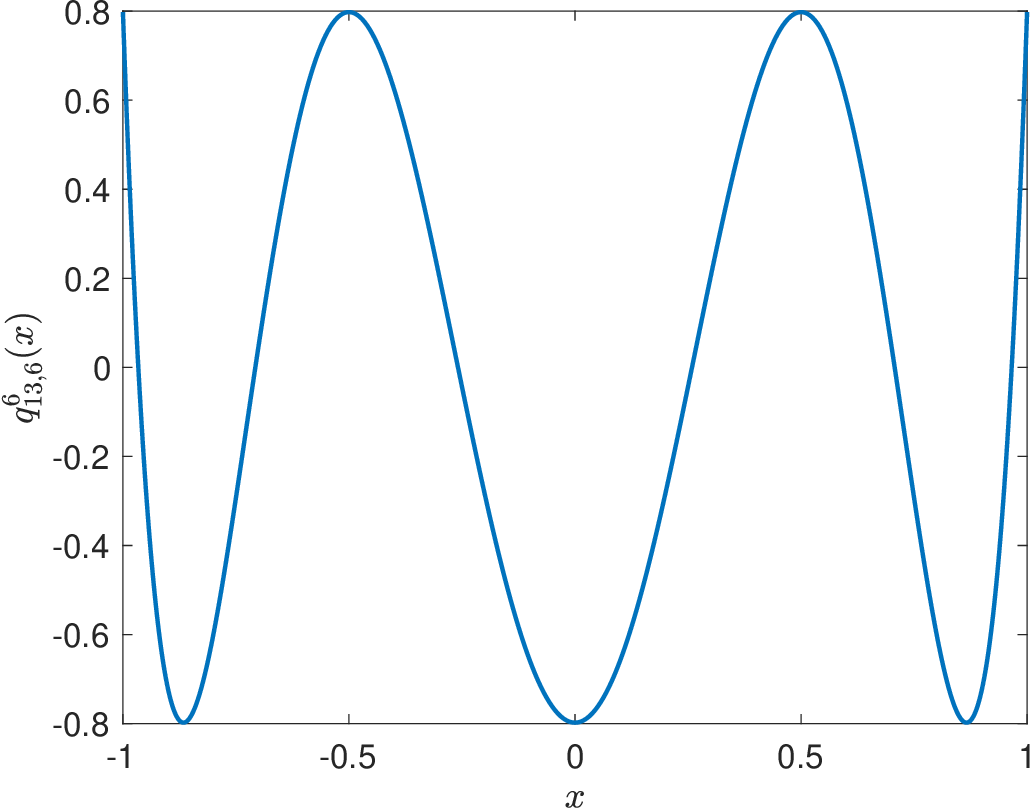}
\hspace{0.2cm}
\includegraphics[scale = 0.40]{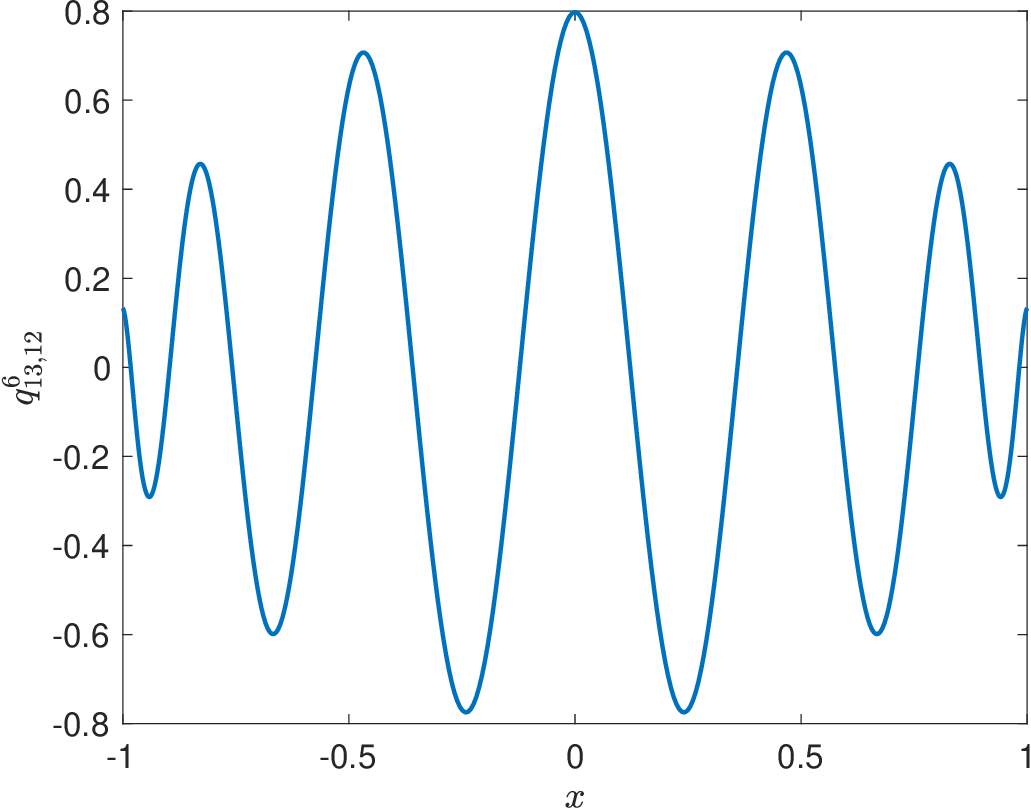}
\vspace{-.5cm}
\end{center}
\caption{ The polynomials $q_{13,6}^6$ (left) and $q_{13,12}^6$  (right)
.\label{fig_q} }
\end{figure}

Note that, like the Chebyshev polynomials, also the orthogonal basis \eqref{q-basis} of $\V_n^m$  is not localized in $[-1,1]$, as shown in Figure~\ref{fig_q} where the polynomials $q_{13,6}^6$ and $q_{13,12}^6$ are plotted.
\begin{definition}\label{def-sca}
Under the previous notation,  we define the orthonormal VP scaling functions as follows
\begin{equation}\label{sca-ort}
\tilde\varphi_{n,k}^m(x):=\sum_{r=0}^{n-1}\tau_{r,k}^{n} q_{n,r}^m(x),\qquad k=1,\ldots,n
\end{equation}
where
\begin{equation}\label{tau}
\tau_{r,k}^{n}:=\sqrt{\frac\pi {n\ \nu_{n,r}^m}}\ p_r(x_k^n), \qquad r=0,\ldots, n-1,\qquad k=1,\ldots, n.
\end{equation}
\end{definition}
\begin{figure}[!htb]%
\begin{center}
\includegraphics[scale = 0.40]{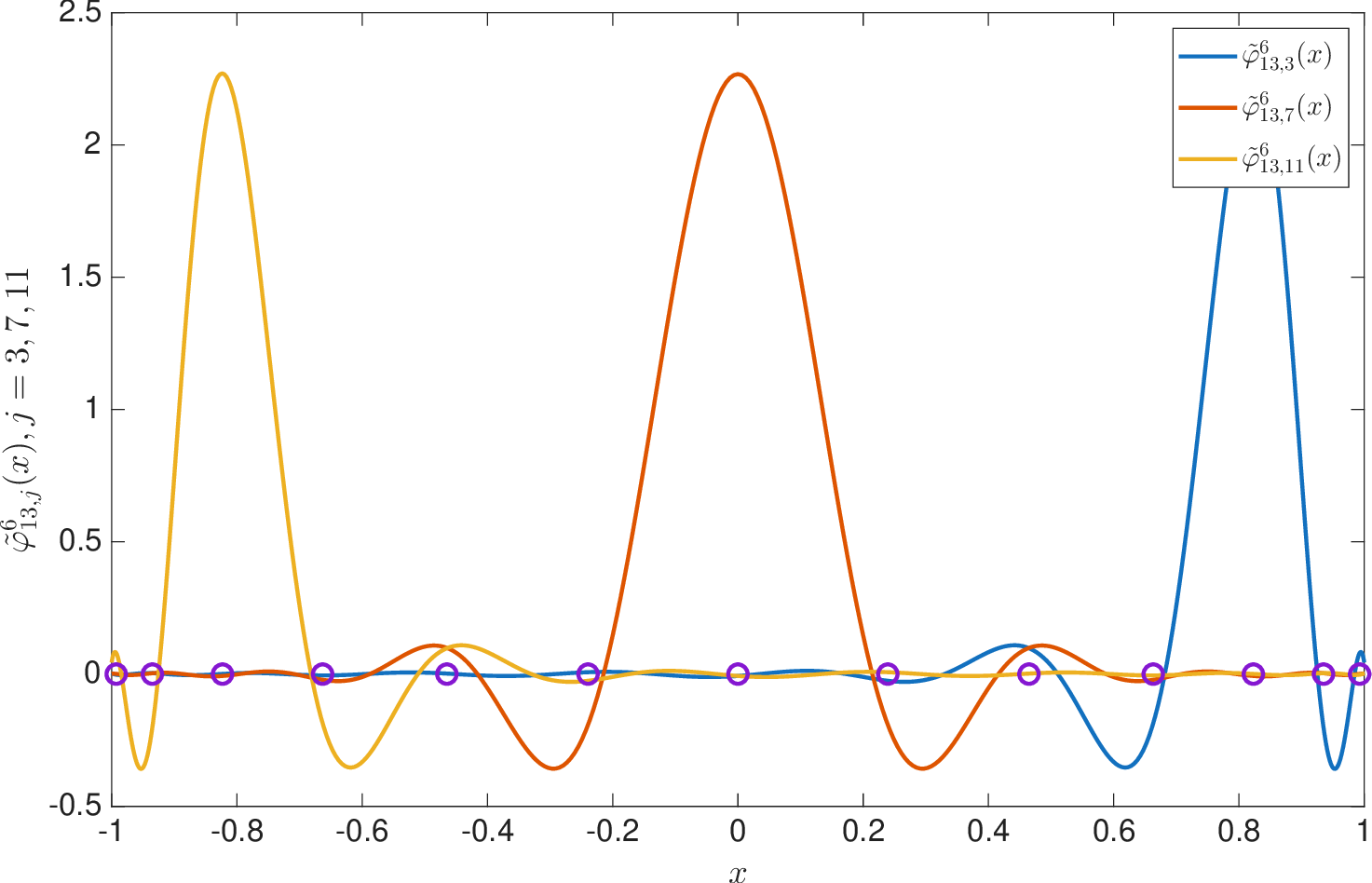}
\vspace{-.5cm}
\end{center}
\caption{ The orthonormal VP scaling functions $\tilde\varphi_{n,k}^m$ for $n=13$, $m=6$, and $k=3,7,11$. The open dots denote the set of nodes $X_n$.\label{fig_phit} }
\end{figure}
Figure~\ref{fig_phit} displays some plots of the VP scaling polynomials previously  defined.
We can observe that, contrary to the orthogonal polynomials in \eqref{q-basis}, the VP scaling polynomials defined in \eqref{sca-ort} exhibit a good localization around the Chebyshev nodes. 
Moreover, they behave very similarly to the interpolatory VP scaling basis defined in \cite{TB-wave} (compare Figures \ref{fig_Phi} and \ref{fig_phit})  but they are not interpolating anymore. In the following proposition, we state that they constitute an orthonormal basis of $\V_n^m$.
\begin{proposition}\label{prop-sca}
We have that $\V_n^m = \mbox{span} \{\tilde\varphi_{n,k}^m, \ k=1,\ldots,n\}$ with
\begin{equation}\label{prod-sca}
<\tilde\varphi_{n,k}^m,\ \tilde\varphi_{n,h}^m>_{L^2_w}=\delta_{k,h},  \qquad k,h=1,2,\ldots, n.
\end{equation}
\end{proposition}

Finally, the next proposition provides the formula for the change of basis from interpolating to orthonormal VP scaling polynomials.
\begin{proposition}\label{prop-sca1}
With respect to the interpolating basis \eqref{sca}, the orthonormal scaling functions can be expanded as follows
\begin{equation}\label{sca-ort1}
\tilde\varphi_{n,k}^m(x)=\sqrt{\frac\pi n}\sum_{h=1}^n\left[\sum_{r=0}^{n-1}\frac{p_r(x_k^n)p_r(x_h^n)}{\sqrt{\nu_{n,r}^m}}\right]\Phi_{n,h}^m(x) \qquad k=1,2,\ldots, n.
\end{equation}
\end{proposition}
\begin{remark}
We note that for the opposite change of basis, as an inverse formula of \eqref{sca-ort1}, we have
\begin{equation}\label{sca-ort1-inv}
\Phi_{n,h}^m(x)=\left(\sqrt{\frac\pi n}\right)^3 \sum_{k=1}^n\left[\sum_{r=0}^{n-1}\sqrt{\nu_{n,r}^m} p_r(x_k^n)p_r(x_h^n)\right]\tilde\varphi_{n,k}^m(x), \qquad h=1,2,\ldots, n.
\end{equation}
This formula can be easily proved by computing the Fourier coefficients $<\Phi_{n,h}^m, \tilde\varphi_{n,k}^m>_{L^2_w} $ by means of \eqref{change-sca}, \eqref{sca-ort}, and \eqref{q-prod}.
\end{remark}
\section{Approximation operators on the sampling space \texorpdfstring{$\V_n^m$}{Vnm}}\label{sec-approx-oper}
In this section,  using the 
new orthonormal VP scaling basis given in Def.~\ref{def-sca}, we introduce the corresponding Fourier projection on $\V_n^m$, denoted by $S_n^m$, and a discrete approximation operator on $\V_n^m$, denoted by $\tilde S_n^m$,  which we deduce from $S_n^m$ by approximating the Fourier-like coefficients.

We study their behavior w.r.t. the following norm
\begin{equation}\label{def-norm}
\|f\|_{p}:=\left\{\begin{array}{cc}
 \displaystyle\sup_{|x|\le 1}|f(x)|    &  p=\infty\\
\displaystyle\left(\int_{-1}^1|f(x)|^pw(x)dx \right)^\frac 1p    & 
1\le p<\infty
\end{array}\right.
\end{equation}
in the functional spaces $L^p_w=\{f \ : \ \|f\|_{p}<\infty\}$, with $1\le p\le \infty$. In the case $p=\infty$ we also consider the space $C^0$ of continuous functions on $[-1,1]$ equipped with the sup norm.

Taking into account that $\PP_{n-m}\subseteq \V_n^m\subseteq\PP_{n+m-1}$, we will compare the approximation error with the error of the best polynomial approximation
\begin{equation}\label{E-best}
E_n(f)_p:=\inf_{P\in\PP_n}\|f-P\|_{p}, \qquad 1\le p\le\infty, \qquad n\in\NN, 
\end{equation}
looking for having a comparable approximation order. 
With regard to this, we recall that 
\begin{equation}\label{lim-En}
  \lim_{n\to\infty}E_n(f)_p=0, \qquad 1\le p\le\infty  
\end{equation}
holds $\forall f\in L^p_w$ in the case $1\le p<\infty$ while it holds $\forall f\in C^0$ if $p=\infty$.\newline  Moreover, the convergence order in \eqref{lim-En} is strictly connected with the smoothness of the function $f$ in the space considered, and it is known for several classes of smoothness \cite{DT}.

That said, as a first operator, let us consider the Fourier projector associated with the new orthonormal VP scaling basis  \eqref{sca-ort} of $\V_n^m$.
As is well known, such an orthogonal projection is defined by
\begin{equation}\label{Sn}
S_n^mf(x):=\sum_{k=1}^n <f,\tilde\varphi_{n,k}^m>_{L^2_w} \tilde\varphi_{n,k}^m(x), \qquad |x|\le 1.
\end{equation}
Equivalently, in integral form, we have
\begin{equation}\label{Sn-int}
S_n^mf(x)=\int_{-1}^1 s_n^m(x,y)f(y)w(y)dy, \qquad s_n^m(x,y):=\sum_{k=1}^n \tilde\varphi_{n,k}^m (x) \tilde\varphi_{n,k}^m(y).
\end{equation}
\begin{remark}
Note that by Definition \ref{def-sca}, using the orthogonality of the matrix $[\sqrt{\frac\pi n}p_{i-1}(x_j^n)]_{i,j=1,\ldots,n}$, we get
\begin{eqnarray*}
s_n^m(x,y)&=&\frac\pi n \sum_{k=1}^n\left(\sum_{r=0}^{n-1}\frac{p_r(x_k^n) q_{n,r}^m(x)}{\sqrt{\nu_{n,r}^m}}\right)
\left(\sum_{s=0}^{n-1}\frac{p_s(x_k^n) q_{n,s}^m(y)}{\sqrt{\nu_{n,s}^m}}\right)\\
&=& 
 \sum_{r=0}^{n-1}\sum_{s=0}^{n-1}\frac{q_{n,r}^m(x)}{\sqrt{\nu_{n,r}^m}}
\frac{q_{n,s}^m(y)}{\sqrt{\nu_{n,s}^m}} \delta_{r,s},
\end{eqnarray*}
i.e., the kernel $s_n^m(x,y)$ can also be written as follows
\begin{equation}\label{sn-q}
 s_n^m(x,y)=\sum_{r=0}^{n-1}\frac{q_{n,r}^m(x)q_{n,r}^m(y)}{\nu_{n,r}^m} , \qquad \forall x,y\in [-1,1]  .
\end{equation}
\end{remark}
 Let us first study the asymptotic behavior of the norm of the map $S_n^m:L^p_w\to L^p_w$ defined by \eqref{Sn}, i.e., let us examine
\[
\|S_n^m\|_{L^p_w\to L^p_w}=\sup_{f\in L^p_w}\frac{\|S_n^mf\|_p}{\|f\|_p}, \qquad 1\le p\le\infty.
\]
The case $p=2$ is trivial since any orthogonal projection in the associated Hilbert space has the norm equal to one, and hence we have
\begin{equation}\label{Sn-norm-2}
\|S_n^m\|_{L^2_w\to L^2_w}=1, \qquad \forall n,m\in\NN\quad \mbox{with $m<n$}.    
\end{equation}
For the other cases, a crucial role is played by the Lebesgue functions $\lambda_n^m(x)$ and the Lebesgue constants $\Lambda_n^m$, defined as
\begin{equation}\label{def-LC}
\lambda_n^m(x) := \int_{-1}^1|s_n^m(x,y)|w(y)dy, \qquad \Lambda_n^m:=\sup_{|x|\le 1} \lambda_n(x) .   
\end{equation} 
Indeed, by classical arguments of functional analysis, for  the limit cases $p=1$ and $p=\infty$, we have
\begin{equation}\label{Sn-norm-1}
\|S_n^m\|_{L^p_w\to L^p_w}=\Lambda_n^m, \qquad p\in \{1,\infty\},
\end{equation}
and using the Riesz–Thorin interpolation theorem (see, e.g., \cite[Corollary 2.2]{Gus}) we also get
\begin{equation}\label{Sn-norm-p}
\|S_n^m\|_{L^p_w\to L^p_w}\le \C_p\ \Lambda_n^m, \qquad 1<p<\infty
\end{equation}
with $\C_p>0$ dependent only on $p$.

\begin{figure}[h!t]%
\begin{center}
\includegraphics[scale = 0.40]{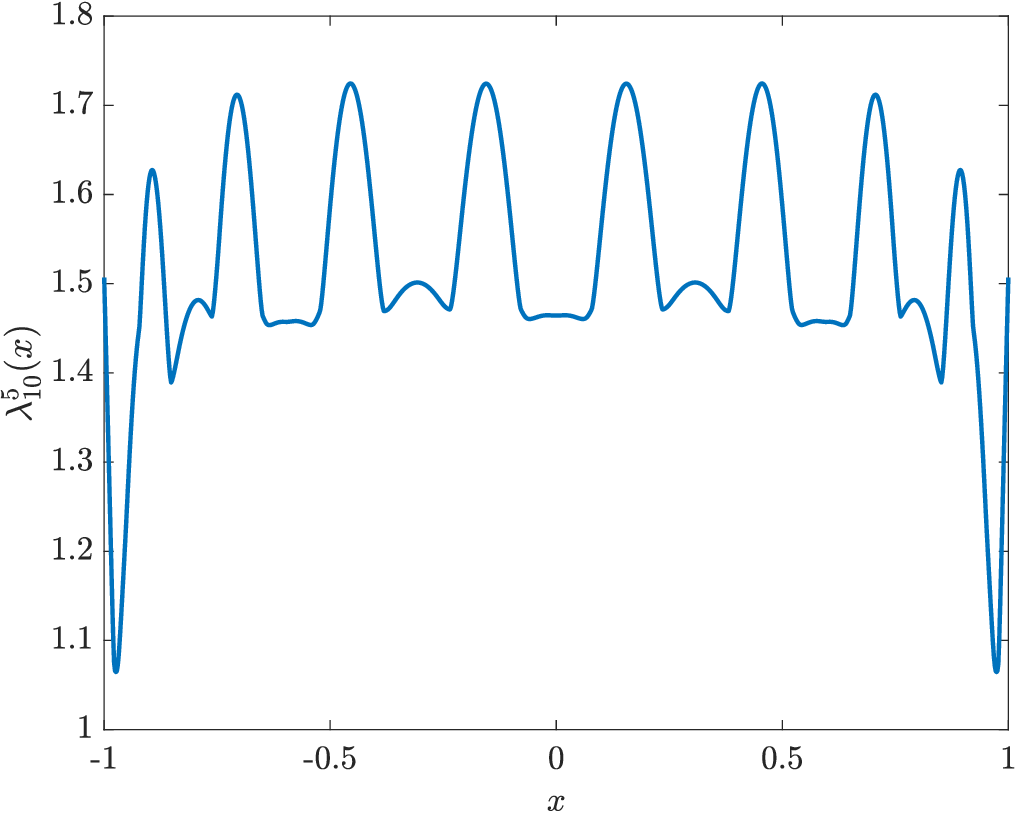}
\hspace{0.2cm}
\includegraphics[scale = 0.40]{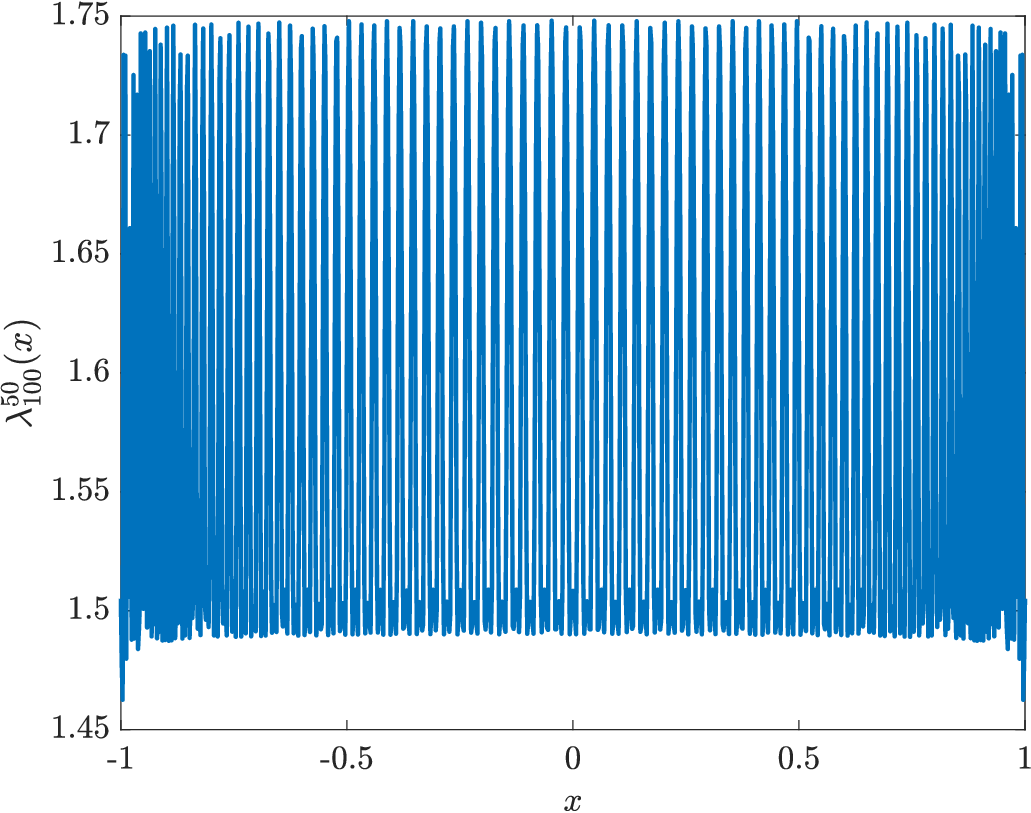}
\vspace{-.5cm}
\end{center}
\caption{ Lebesgue functions $\lambda_{10}^{5}(x)$ (left) and $\lambda_{100}^{50}(x)$ (right). \label{fig_LF_s} }
\end{figure}

In Figure~\ref{fig_LF_s} we show the Lebesgue functions $\lambda_{n}^{m}(x)$
for $n=10,100$ and $m=\lfloor \theta n\rfloor$ with  $\theta=0.5$.

The following theorem states that, for arbitrarily large $n\in\NN$, the Lebesgue functions $\lambda_n^m(x)$ are bounded in $[-1,1]$, uniformly w.r.t. $n$.
\begin{theorem}\label{th-LC}
Let $\theta\in ]0,1[$ arbitrarily fixed, $n\in\NN$ arbitrarily large, and $m=\lfloor \theta n\rfloor$.  The Lebesgue constants defined in \eqref{def-LC} satisfy
\begin{equation}\label{Sn-LC}
\sup_n \Lambda_n^m <\infty .
\end{equation}
\end{theorem}

As a corollary of the previous result, we have the following
\begin{corollary}\label{cor-Snm}
Let $\theta\in ]0,1[$ arbitrarily fixed, $n\in\NN$ arbitrarily large, and $m=\lfloor \theta n\rfloor$. For any $1\le p\le \infty$, the  operator map $S_n^m:L^p_w\to L^p_w$  defined by \eqref{Sn} is uniformly bounded w.r.t. $n$. Moreover,  we have
\begin{equation}\label{Snm-nearbest}
 E_{n+m-1}(f)_p\le \|S_n^mf-f\|_p\le \C E_{n-m}(f)_p, \qquad \forall f\in L^p_w,\qquad 1\le p\le \infty   
\end{equation}
where $\C>0$ is a constant independent of $f,n,m$.
\end{corollary}
\begin{remark}\label{rem-Sn}
Regarding the mapping properties of $S_n^m:L^p_w\to L^p_w$, we point out that, by the previous arguments, for an arbitrary fixed pair of positive integers $n>m$, and for $1\le p\le\infty$, we have
 \begin{equation}\label{eq-remS}
     \|S_n^mf\|_p\le \C \|f\|_p\qquad \mbox{with}\qquad \C=\left\{\begin{array}{ll}
         \Lambda_n^m &  p=1,\infty\\
          1 & p=2\\
          \C_p\Lambda_n^m & \mbox{otherwise}.
     \end{array}\right.
 \end{equation}
 being $\C_p>0$ dependent only on $p$.
\newline 
Moreover, taking into account that
\[
\|S_n^mf- S_n^m\tilde f \|_p\le  \|S_n^m\|_{L^p_w\to L^p_w} \|f-\tilde f\|_p, \qquad 1\le p\le \infty
\]
the uniform boundedness of the operator $S_n^m: L^p_w\to L^p_w$ ensures the well--conditioning. 
\newline
Finally, regarding the consequences of the error estimate \eqref{Snm-nearbest}, we note that, as $n\to \infty$, the choice  $m=\lfloor \theta n\rfloor$ with fixed $\theta\in ]0,1[$, ensures that both $(n-m)$ and $(n+m-1)$ tend to $\infty$ at the same rate of $n$. Consequently, the estimate \eqref{Snm-nearbest} implies that
\begin{equation}\label{conv-order-S}
 \|S_n^mf-f\|_p =\bigO(n^{-a}) \Longleftrightarrow E_n(f)_p=\bigO(n^{-a}), \qquad \forall a>0, \qquad 1\le p\le\infty   .
\end{equation}
\end{remark}

\begin{figure}[!ht]%
\begin{center}
\includegraphics[scale = 0.40]{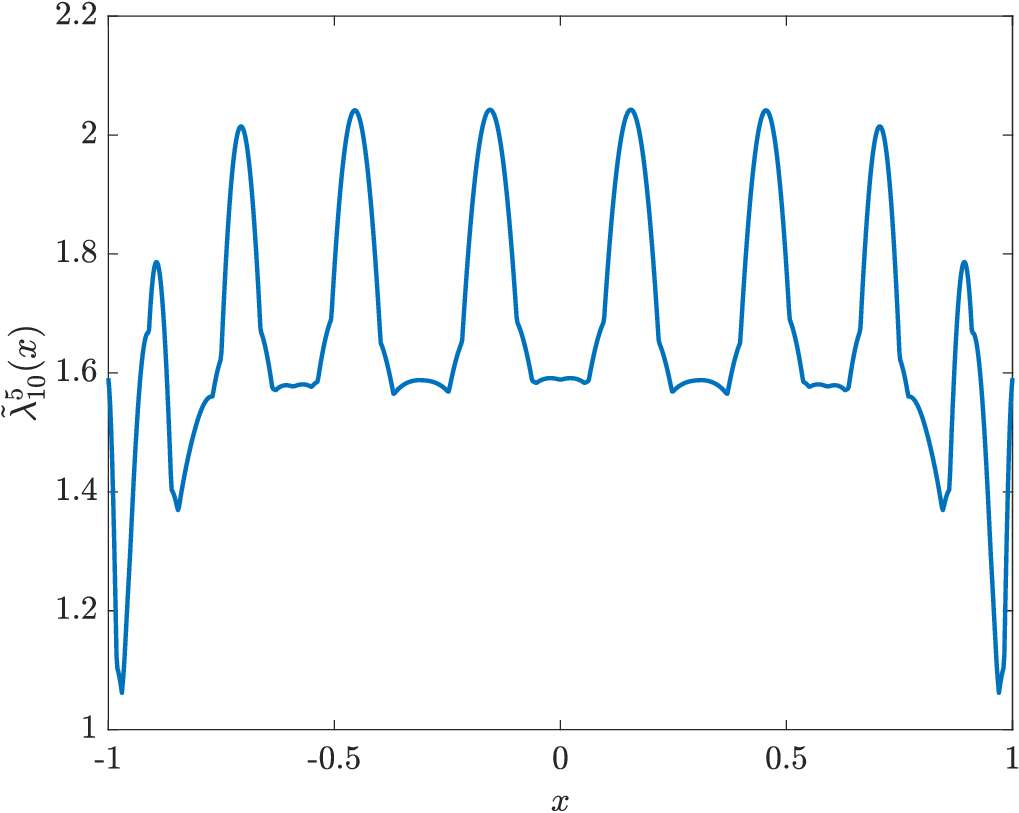}
\hspace{0.2cm}
\includegraphics[scale = 0.40]{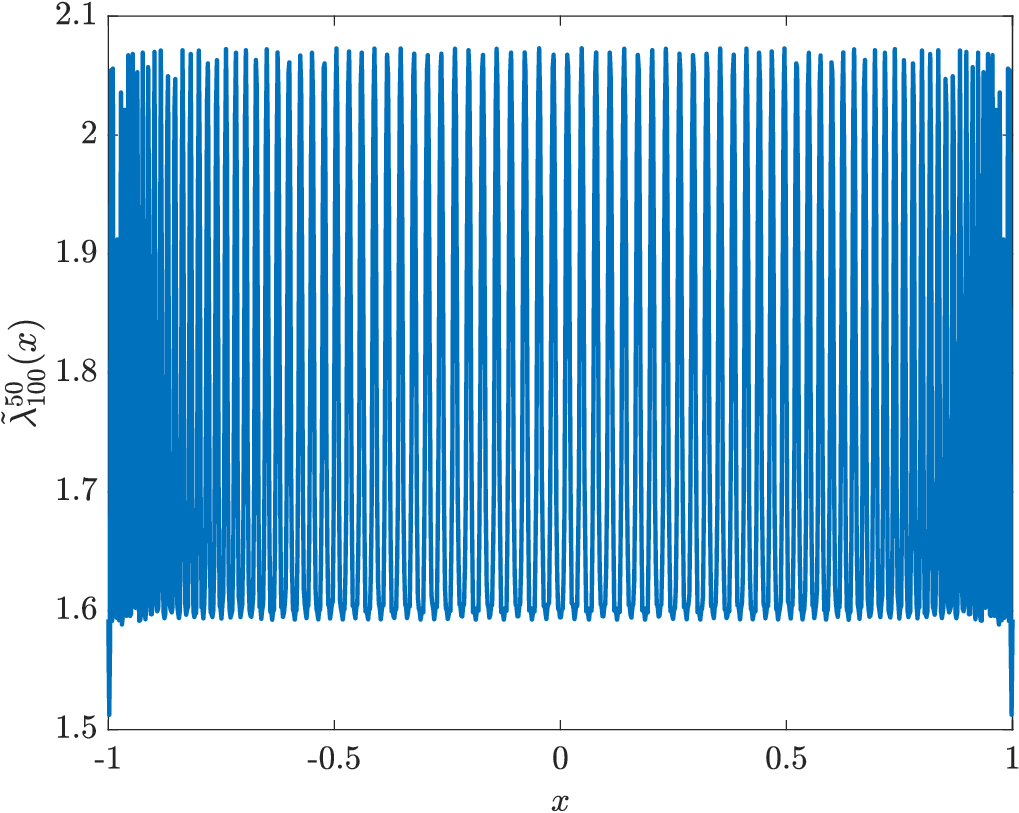}
\vspace{-.5cm}
\end{center}
\caption{ Lebesgue functions $\tilde\lambda_{10}^{5}(x)$ (left) and $\tilde\lambda_{100}^{50}(x)$ (right)
.\label{fig_LF_st} }
\end{figure}

By the previous results the Fourier projection $S_n^mf$ provides a near best and well--conditioned approximation of any $f\in L^p_w$ in the polynomial space $\V_n^m$. However, to obtain $S_n^mf$ we need to compute the integral coefficients in \eqref{Sn}, which can be a difficult task. 
To overcome this problem,  we introduce a discrete version of the orthogonal projection $S_n^mf$, achieved by applying to the Fourier--like coefficients in \eqref{Sn} the Gauss-Chebyshev quadrature rule on $n$ nodes
\begin{equation}\label{Gauss}
 \int_{-1}^1 P(x) w(x) dx =\frac \pi n \sum_{i=1}^n  P(x_i^n), \qquad \forall P\in\PP_{2n-1}  .
\end{equation}
In this way, we define the following discrete version of $S_n^mf$
\begin{equation}\label{Sn-tilde}
\tilde S_n^mf(x):=\sum_{k=1}^n \left[\frac\pi n \sum_{i=1}^n f(x_i^n)\tilde\varphi_{n,k}^m(x_i^n)\right] \tilde\varphi_{n,k}^m(x), \qquad |x|\le 1.
\end{equation}
Equivalent forms of this polynomial are the following
\begin{eqnarray}\label{Sn-tilde1}
\tilde S_n^mf(x)&=&\frac\pi n \sum_{i=1}^n f(x_i^n)s_n^m(x_i^n, x),\\
\label{Sn-tilde2}
\tilde S_n^mf(x)&=&\sum_{k=1}^n\left[\left(\frac\pi n\right )^\frac 32\sum_{h=1}^n f(x_h^n) \sum_{r=0}^{n-1}\frac{p_r(x_k^n)p_r(x_h^n)}{\sqrt{\nu_{n,r}^m}}\right]\tilde\varphi_{n,k}^m(x).
\end{eqnarray}
For the discrete projection $\tilde S_n^m$ defined by the previous identities,  the Lebesgue functions and constants are the following
\begin{equation}\label{def-LC1}
 \tilde\lambda_n^m(x) := \frac\pi n \sum_{i=1}^n |s_n^m(x_i^n, x)|, \qquad
 \tilde\Lambda_n^m := \sup_{|x| \leq 1} \tilde\lambda_n^m(x).
\end{equation}
Similarly to the case of the Fourier projector, we have 
\begin{equation}\label{tildeSn-LC-inf}
\|\tilde S_n^m\|_{C^0\to C^0}= \tilde\Lambda_n^m, \qquad \forall n>m 
\end{equation}
and, more specifically, by the definition, for all functions $f$  we get
\begin{equation}\label{LC-point}
|\tilde S_n^mf(x)|\le \tilde \lambda_n (x)
\left(\max_{1\le i\le n}|f(x_i^n)|\right), \qquad |x|\le 1.
\end{equation}
Figure~\ref{fig_LF_st} shows  the Lebesgue functions
$\tilde\lambda_{n}^{m}(x)$ for $n=10,100$ and $m=n/2$. 

\begin{figure}[!htb]%
\begin{center}
\includegraphics[scale = 0.40]{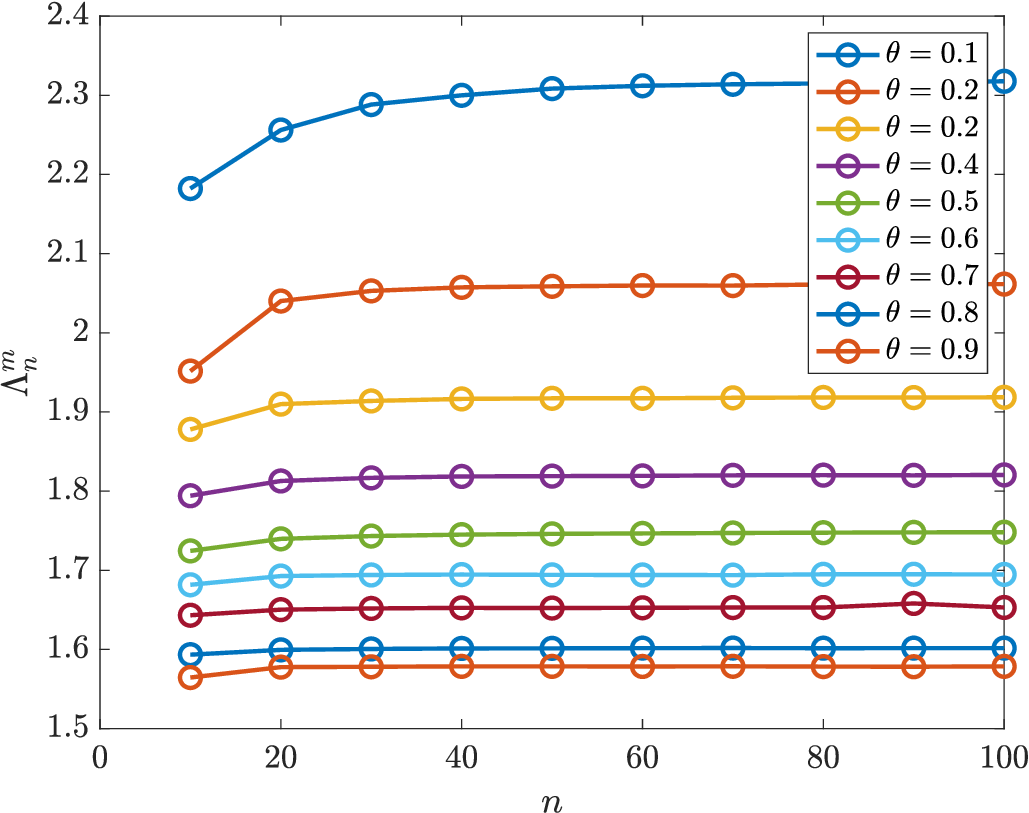}
\hspace{0.2cm}
\includegraphics[scale = 0.40]{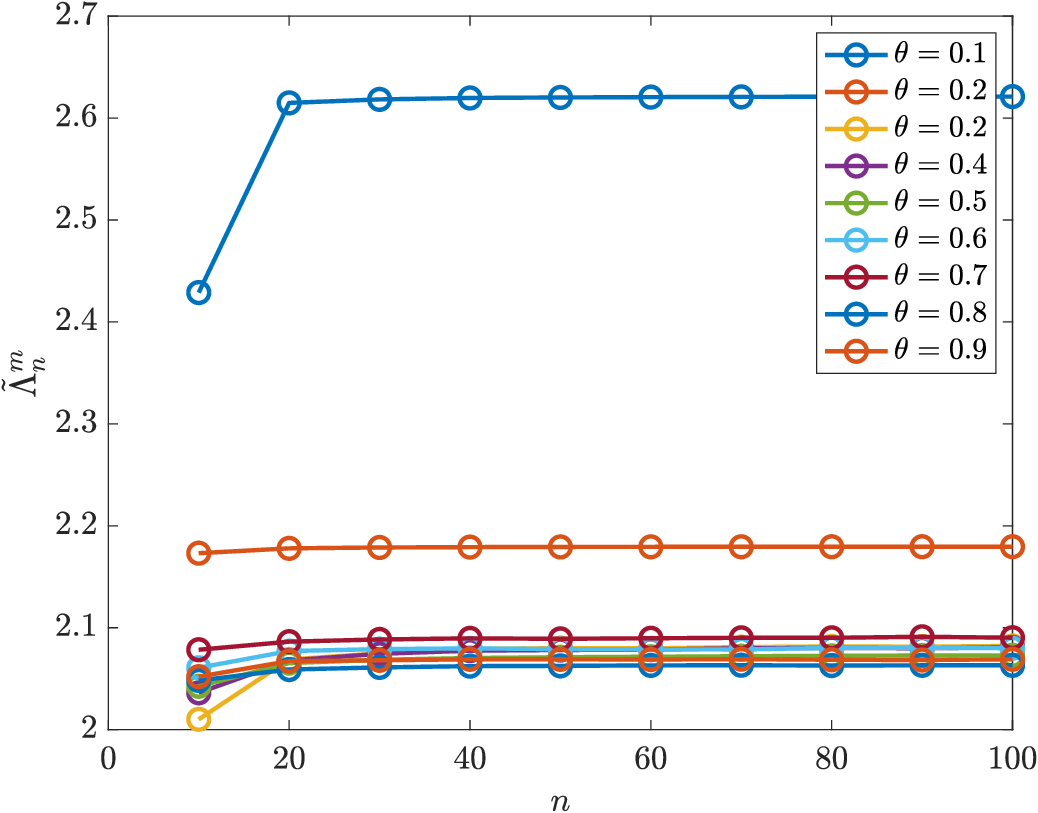}
\vspace{-.5cm}
\end{center}
\caption{ Lebesgue constants $\Lambda_n^m$ (left) and
$\tilde\Lambda_m^n$ (right) for $n = 10,20,\ldots,100$ and $m=\lfloor \theta n\rfloor$ with $\theta = 0.1,0.2,\ldots,0.9$.
\label{fig_LC} }
\end{figure}

Comparing Figure~\ref{fig_LF_st} with Figure \ref{fig_LF_s} we note the similar behavior of the Lebesgue functions related to $\tilde S_n^m$ and $S_n^m$ when $n=10,100$ and $m=n/2$.
Indeed, such a similar behavior holds more generally for all $n\in\NN$, as stated below.
\begin{theorem}\label{th-LF}
For all $n,m\in\NN$ such that $m=\lfloor \theta n\rfloor$ with fixed $\theta\in ]0,1[$, the Lebesgue functions defined in \eqref{def-LC} and \eqref{def-LC1} satisfy
\begin{equation}\label{eq-LF}
   \C_1 \lambda_n^m(x)\le \tilde\lambda_n^m(x)\le \C_2 \lambda_n^m(x), \qquad |x|\le 1,
\end{equation}
where $\C_1,\C_2>0$ are constants independent of $n,m,x$.
\end{theorem}
Note that, taking the supremum w.r.t. $x\in [-1,1]$ in \eqref{eq-LF}, we obtain the following inequalities for the associated Lebesgue constants
\begin{equation}\label{eq-LF1}
    \C_1\Lambda_n^m\le \tilde\Lambda_n^m\le \C_2\Lambda_n^m.
\end{equation}
 Figure~\ref{fig_LC} shows the Lebesgue constants $\Lambda_n^m$ and
$\tilde\Lambda_m^n$ for $n = 10,20,\ldots,100$ and $m=\lfloor \theta n\rfloor$, with $\theta = 0.1,0.2,\ldots,0.9$. Relation \eqref{eq-LF1} is confirmed with $\C_1=1$ and we conjecture that $\Lambda_n^m\le \tilde \Lambda_n^m$ is always true.

Concerning the behavior of the operator $\tilde S_n^m$ in spaces $L^p_w$ with $1\le p\le \infty$, we have to take into account the discrete nature of $\tilde S_n^m$ and, in addition to considering functions defined everywhere, we must replace the continuous norm $\|f\|_p$ with the following discrete norm
\[
\|f\|_{\ell^p(X_n)}:=\left\{\begin{array}{cc}
\displaystyle \max_{1\le k\le n}|f(x_k^n)|    &  p=\infty\\
\displaystyle\left(\frac \pi n \sum_{k=1}^n|f(x_k^n)|^p \right)^\frac 1p    & 
1\le p<\infty
\end{array}\right.\qquad n\in\NN
\]
in order to get some uniform boundedness result.

We point out that, in general, $\|f\|_{\ell^p(X_n)}$ is  not a norm, but restricted to the polynomials $f\in \V_n^m$ it is a norm since by \eqref{sca-inter}-\eqref{Sca-space}, $\forall f,g\in\V_n^m$, we have
\begin{equation}\label{id}
 f(x)=g(x), \quad \forall x\in [-1,1]\ \Longleftrightarrow \ f(x_i^n)=g(x_i^n),\quad i=1,\ldots,n,    
\end{equation}
which ensures that $\|f\|_{\ell^p(X_n)}=0 \Longrightarrow f=0$,  $\forall f\in \V_n^m$.

Using the previous discrete norms, we get the following analog of Corollary \ref{cor-Snm}.
\begin{theorem}\label{th-tildeSn}
 Let $\theta\in ]0,1[$ arbitrarily fixed, $n\in\NN$ arbitrarily large, and $m=\lfloor \theta n\rfloor$.  We have
\begin{equation}\label{tildeSn-p}
\|\tilde S_n^m f\|_p\le\C_p \|f\|_{\ell^p(X_n)}, \qquad  1\le p\le \infty,
\end{equation}
where $\C_p>0$ is a constant that depends only on $p$.
\newline
Moreover, for all $f\in C^0$, the polynomial $\tilde S_n^mf$ uniformly converges to $f$ as $n\to \infty$, and the  near-best error estimate
\begin{equation}\label{Stilde-nearbest}
 E_{n+m-1}(f)_\infty\le \|\tilde S_n^m f-f\|_\infty\le \C E_{n-m}(f)_\infty   
\end{equation}
holds with $\C>0$ independent of $f,n,m$.
\end{theorem}
\begin{remark}\label{rem-tildeSn}
Regarding the value of the constant $\C_p$ in \eqref{tildeSn-p}, in proving Theorem \ref{th-tildeSn} we state that \eqref{eq-remStilde} holds for any pair of integers $0<m<n$.

Finally, recalling the discussion of Remark  \ref{rem-Sn}, we observe that the error estimate \eqref{Stilde-nearbest} implies the following equivalence
\begin{equation}\label{conv-order-Stilde}
 \|\tilde S_n^mf-f\|_\infty =\bigO(n^{-a}) \Longleftrightarrow E_n(f)_\infty=\bigO(n^{-a}), \qquad \forall a>0,  
\end{equation}
holds as $n\to\infty$ and $m=\lfloor \theta n\rfloor$ with fixed $\theta\in ]0,1[$.
\end{remark}
In conclusion, we have stated that the polynomial $\tilde S_n^mf$, similarly to the orthogonal projection $S_n^mf$, provides a well--conditioned, near-best uniform approximation in the sampling space $\V_n^m$, and, contrary to $S_n^mf$, it is easily computable once the values of $f$  are known at the Chebyshev nodes $X_n$. 
\section{Orthonormal VP wavelet functions}\label{sec-wavelets}
\begin{figure}[!b]%
\begin{center}
\includegraphics[scale = 0.40]{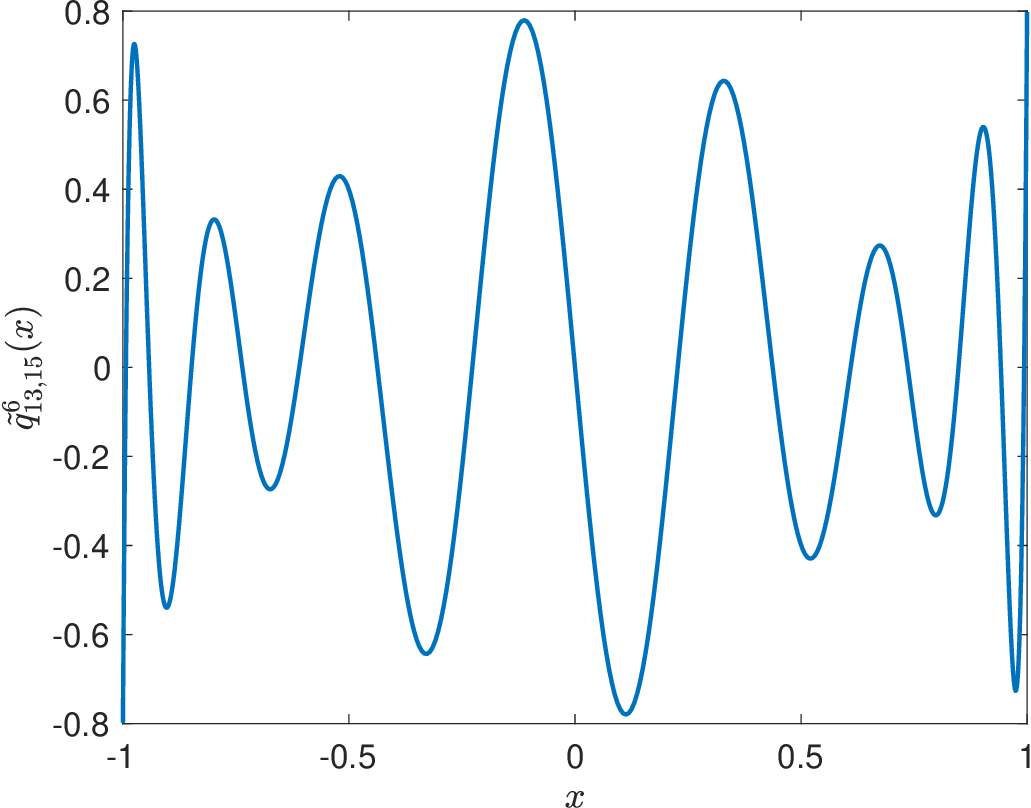}
\hspace{0.2cm}
\includegraphics[scale = 0.40]{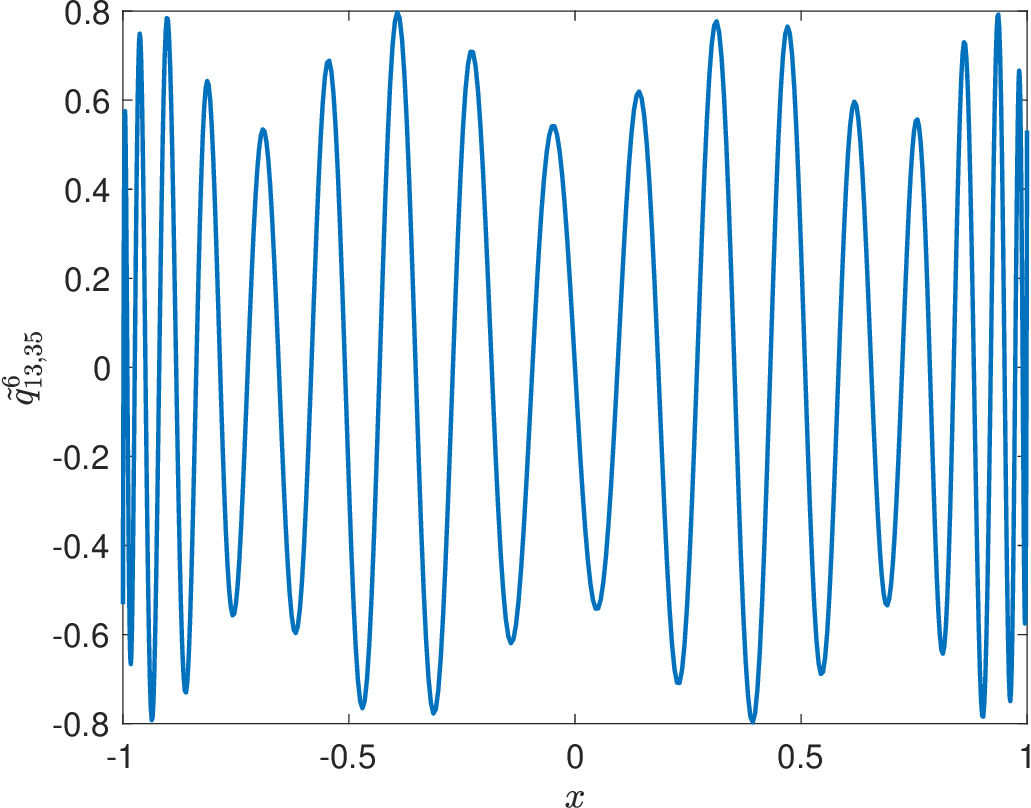}
\vspace{-.7cm}
\end{center}
\caption{ The polynomials $\tilde{q}_{13,15}^6$ (left) and $\tilde{q}_{13,35}^6$ (right).
.\label{fig_qt} }
\end{figure}
Firstly, we recall that in the wavelet space $\W_n^m$
we have already found the orthogonal basis \cite{TB-wave}
\begin{equation}\label{wav-ort}
\tilde q_{n,r}^m(x):=\left\{\begin{array}
{ll}
\displaystyle \mu_{n, r}^m p_{2n-r}(x)+ \mu_{n, 2n-r}^m p_r(x) & n\le r<n+m\\ [.1in]
p_r(x)& n+m\le r\le 3n-m\\ [.1in]
\displaystyle q_{3n,r}^m(x)= \mu_{3n, r}^m p_{r}(x)- \mu_{3n, 6n-r}^m p_{6n-r}(x) & 3n-m< r<3n
\end{array}\right.
\end{equation}
satisfying 
\begin{equation}\label{prod1}
<\tilde q_{n,r}^m, \tilde q_{n,s}^m>_{L_w^2}=v_{n,r}^m \delta_{r,s},\qquad r,s=n,\ldots, 3n-1 ,
\end{equation}
\begin{equation}\label{vr}
 v_{n,r}^m:=\left\{\begin{array}
{ll}
\displaystyle \frac{m^2+(n-r)^2}{2m^2} & n< r<n+m\\ [.1in]
1& r=n \quad\mbox {or}\quad n+m\le r\le 3n-m\\ [.1in]
\displaystyle \frac{m^2+(3n-r)^2}{2m^2} & 3n-m< r<3n.
\end{array}\right. 
\end{equation}
In Figure~\ref{fig_qt} we show the plots of the polynomials $\tilde{q}_{13,15}^6$ and $\tilde{q}_{13,35}^6$. As we can see, the orthogonal basis polynomials \eqref{wav-ort} are not localized.

\begin{figure}[!hb]%
\begin{center}
\includegraphics[scale = 0.40]{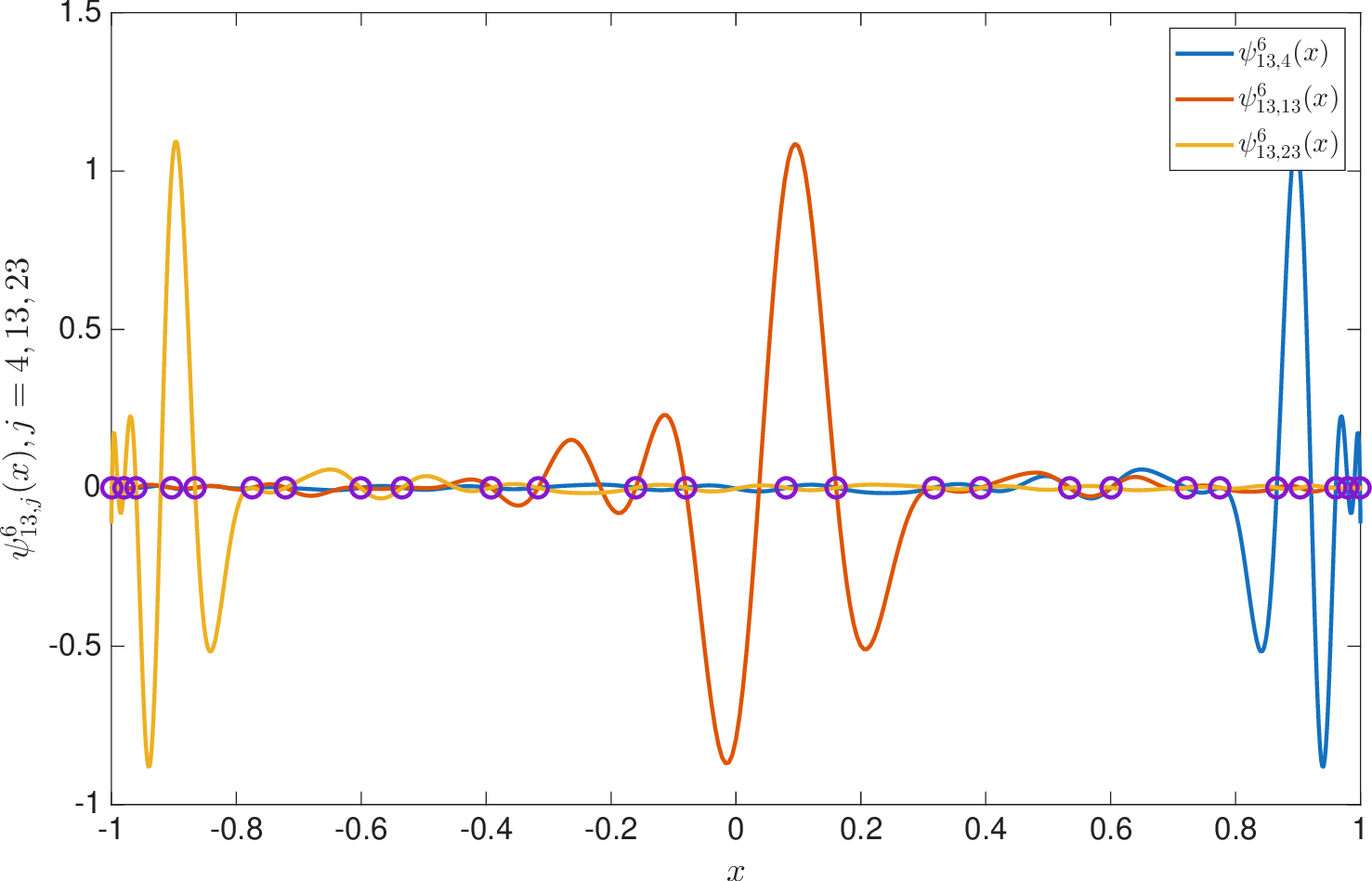}
\vspace{-.5cm}
\end{center}
\caption{The interpolating VP wavelets $\psi_{n,k}^m$ for $n=13$, $m=6$, and $k=4,13,23$. The open dots denote the set of nodes $Y_{2n}$.\label{fig_psi} }
\end{figure}

However, in \cite{TB-wave} the authors introduce well localized polynomial wavelets $\psi_{n,k}^m$ that for $k=1,\ldots, 2n$ constitute a Riesz basis of the detail space $\W_n^m$. These wavelets are not orthogonal, but they are interpolating at the nodes
\[
Y_{2n}=X_{3n}-X_n=\{y_k^n, \ k=1,\ldots,2n\}, \qquad n\in\NN,
\]
namely \eqref{wave-inter} holds. 
 They can be expanded in the orthogonal basis \eqref{wav-ort} as follows \cite{TB-wave}
\begin{equation}\label{bchange-wav1}
\psi_{n,k}^m(x)=\frac\pi{3n}\sum_{r=n}^{3n-1}\rho_{r,k} \ \tilde q_{n,r}^m (x),
\qquad k=1,\ldots, 2n
\end{equation}
where
\begin{equation}\label{bchange-wav2}
\rho_{r,k}:=\left\{\begin{array}{ll}
p_n(y_k^n) & r=n\\ [.1in]
p_r(y_k^n)+p_{|2n-r|}(y_k^n) & n<r\le 3n-m \quad r\ne 2n\\ [.1in]
p_{2n}(y_k^n)+\sqrt{2}p_0(y_k^n) & r=2n\\ [.1in]
p_r(y_k^n)+\mu_{n,r-2n}^m p_{r-2n}(y_k^n)-\mu_{n,4n-r}^m p_{4n-r}(y_k^n) & 3n-m<r<3n .
\end{array}\right.
\end{equation}
In Figure~\ref{fig_psi} we show the plots of the VP wavelets $\psi_{13,4}^6$,
$\psi_{13,13}^6$ and $\psi_{13,23}^6$, interpolating at the set of nodes $Y_{26}$ marked with circles.

\begin{figure}[!t]%
\begin{center}
\includegraphics[scale = 0.40]{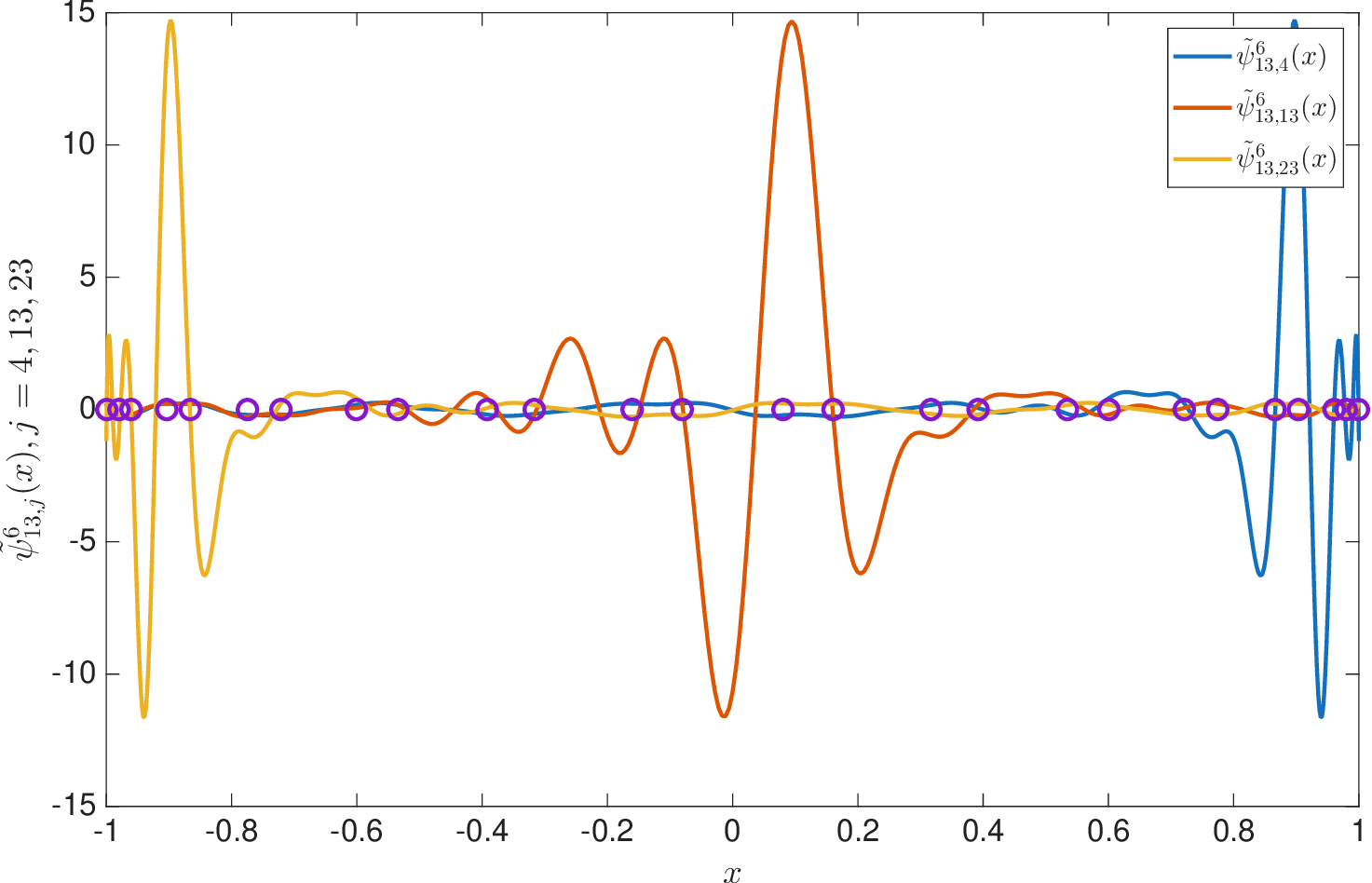}
\vspace{-.5cm}
\end{center}
\caption{The orthonormal VP wavelets $\tilde\psi_{n,k}^m$ for $n=13$, $m=6$, and $k=4,13,23$. The open dots denote the set of nodes $Y_{2n}$.\label{fig_psit} }
\end{figure}

In the following, we define a new wavelet basis of the detail space $\W_n^m$ that is both well-localized and orthonormal.
\begin{definition}\label{def-wav}
Under the previous notation,  we define the orthonormal VP wavelet functions as follows
\begin{equation}\label{wav-ort2}
\tilde\psi_{n,k}^m(x):=\sum_{r=n}^{3n-1}\sigma_{r,k}^{n}\ \frac{\tilde q_{n,r}^m(x)}{\sqrt{v_{n,r}^m}},\qquad k=1,\ldots,2n
\end{equation}
\begin{equation}\label{sigma}
\sigma_{r,k}^{n}:=\sqrt{\frac\pi {3n}}\left\{
\begin{array}{ll}
p_n(y_k^n) & \mbox{if $r=n$}\\ [.15in]
\displaystyle \frac{p_r(y_k^n)+p_{2n-r}(y_k^n)}{\sqrt{2}} & \mbox{if $n<r<2n$}\\[.15in]
\displaystyle \frac{p_{2n}(y_k^n)+\sqrt{2}\ p_0(y_k^n)} {\sqrt{3}} &\mbox{if $r=2n$}\\ [.15in]
\displaystyle \sqrt{\frac 32}\ p_r(y_k^n) &\mbox{if $2n<r<3n$}
\end{array}\right. \qquad k=1,\ldots, 2n.
\end{equation}
\end{definition}
In Figure~\ref{fig_psit} we show the plots of the wavelets $\tilde\psi_{13,4}^6$,
$\tilde\psi_{13,13}^6$ and $\tilde\psi_{13,23}^6$. 
We note similarities with Figure \ref{fig_psi} even if the interpolation property on the set of nodes $Y_{26}$ (marked with empty dots) is no longer valid.
\begin{proposition}\label{prop-wav}
We have that $\W_n^m = \mbox{span} \{\tilde\psi_{n,k}^m, \ k=1,\ldots,2n\}$ with
\begin{equation}\label{prod-wav}
<\tilde\psi_{n,k}^m,\ \tilde\psi_{n,h}^m>_{L^2_w}=\delta_{k,h},  \qquad k,h=1,\ldots, 2n.
\end{equation}
\end{proposition}

\section{Decomposition and reconstruction}\label{sec-dec-and-rec}
As is well-known, the relation \eqref{compl} between the approximation and detail spaces establishes a one-to-one correspondence  between any function $f_{3n}\in \V_{3n}^m$ to which a unique pair of functions $f_n\in\V_n^m$ and $g_{n}\in\W_n^m$ corresponds, such that
\begin{equation}\label{f3n}
    f_{3n}(x)=f_n(x)+g_{n}(x),\qquad \forall x\in [-1,1].
\end{equation}
On the other hand, each function in \eqref{f3n} is uniquely determined by its coefficients in the orthonormal VP scaling or wavelet bases previously defined. To distinguish such coefficients at any resolution level $n\in\NN$, we use the following notation 
\begin{eqnarray}\label{fn}
    f_n(x)&=&\sum_{k=1}^n a_{n,k}\tilde\varphi_{n,k}^m(x), \qquad  f_n\in\V_n^m\\
    \label{gn}
    g_n(x)&=&\sum_{h=1}^{2n} b_{n,k}\tilde\psi_{n,k}^m(x), \qquad  g_n\in\W_n^m
\end{eqnarray}
and we set
\[
\vec{a}_n:=(a_{n,1},\ldots, a_{n,n})^T,\qquad \vec{b}_n:=(b_{n,1},\ldots, b_{n,2n})^T.
\]
The following theorem provides the decomposition/reconstruction formulas to compute $\vec{a}_n$ and $\vec{b}_n$ from $\vec{a}_{3n}$ and vice versa.
\begin{theorem}\label{th-dec}
For each resolution level $n\in\NN$ and any parameter $0<m<n$, the basis coefficients, defined by \eqref{fn}--\eqref{gn}, of the functions 
$f_{3n}=f_n+g_n$ are related by the matrices $A=(A_{kj})\in\RR^{n\times 3n}$ and $B=(B_{hj})\in\RR^{2n\times 3n}$ which are defined by means of the coefficients in \eqref{tau}, \eqref{muj}, \eqref{sigma}, and \eqref{vr}, as follows
\begin{eqnarray}\label{A}
 A_{k,j}&:=& \sum_{s=0}^{n-m} \tau_{s,k}^n \tau_{s,j}^{3n}+
\sum_{s=n-m+1}^{n-1}\tau_{s,k}^n\left[\mu_{n,s}^m \tau_{s,j}^{3n}
-\mu_{n,2n-s}^m \tau_{2n-s,j}^{3n}\right]  
\\
\nonumber
&& \hspace{5.4cm} k=1,\ldots, n,\quad j=1,\ldots,3n
\\ 
\nonumber
&&
\\
\label{B} 
B_{hj}&:=& 
\sum_{r=n}^{n+m-1}\frac{\sigma_{r,h}^n}{\sqrt{v_{n,r}^m}}\ \left[\mu_{n,r}^m \tau_{2n-r,j}^{3n}+\mu_{n,2n-r}^m \tau_{r,j}^{3n}\right]+
\sum_{r=n+m}^{3n-m}\frac{\sigma_{r,h}^n}{\sqrt{v_{n,r}^m}}\ \tau_{r,j}^{3n}
\\
\nonumber
&+&
 \sum_{r=3n-m+1}^{3n-1}\sigma_{r,h}^n\sqrt{v_{n,r}^m} \ \tau_{r,j}^{3n}
 \qquad\qquad h=1,\ldots, 2n,\quad j=1,\ldots,3n.
\end{eqnarray}
These matrices provides the following 
decomposition and reconstruction formulas
\begin{equation}\label{dec-a}
\vec{a}_n=A\ \vec{a}_{3n},\qquad \vec{b}_n=B\  \vec{a}_{3n}
\end{equation}
\begin{equation}\label{rec-a}
\vec{a}_{3n}=A^T\ \vec{a}_n+B^T\ \vec{b}_n.
\end{equation}
\end{theorem}
 Looking at the scalar form of the decomposition/reconstruction formulas \eqref{dec-a}--\eqref{rec-a} we can note they are primarily based on the following four kinds of transformation:
 \begin{itemize}
     \item {\bf Transform $\T_n$:} $\vec{u}=(u_1,\ldots, u_{n})^T\in \RR^n\rightarrow \T_n \vec{u}= \vec{t}=(t_0,\ldots,t_{n-1})^T\in\RR^n$ defined by
     \begin{equation}\label{T}
         t_r:=\sum_{k=1}^n u_k\ \tau_{r,k}^n,\qquad r=0,\ldots, n-1.
     \end{equation}
     \item {\bf Transform $\T_n^{\prime}$:} $\vec{t}=(t_0,\ldots,t_{n-1})^T\in \RR^n\rightarrow \T_N^{\prime} \vec{t}= \vec{u}=(u_1,\ldots, u_{n})^T\in\RR^n$ defined by
     \begin{equation}\label{T-1}
         u_k:=\sum_{r=0}^{N-1} t_r\ \tau_{r,k}^n,\qquad k=1,\ldots, n.
     \end{equation}
     \item {\bf Transform $\Sigma_n$:} $\vec{u}=(u_1,\ldots, u_{2n})^T\in \RR^{2n}\rightarrow \Sigma_n \vec{u}= \vec{s}=(s_n,\ldots,s_{3n-1})^T\in\RR^{2n}$ defined by
     \begin{equation}\label{S}
         s_r:=\sum_{h=1}^{2n} u_h\sigma_{r,h}^n,\qquad r=n,\ldots, 3n-1.
     \end{equation}
     \item {\bf Transform $\Sigma_n^{\prime}$:} $\vec{s}=(s_n,\ldots,s_{3n-1})^T\in \RR^{2n}\rightarrow \Sigma_n^{\prime} \vec{s}= \vec{u}=(u_1,\ldots, u_{2n})^T\in\RR^{2n}$ defined by
     \begin{equation}\label{S-1}
         u_h:=\sum_{r=n}^{3n-1} s_r \sigma_{r,h}^n \qquad h=1,\ldots,2n.
     \end{equation}
 \end{itemize}
  Indeed, the scalar version of the formulas \eqref{dec-a} and \eqref{rec-a} involves double summations, and changing the order of such summations  we easily get the next decomposition and reconstruction algorithms based on the previous transformations and on the parameters  $\mu_r:=\mu_{n,r}^m$ (cf.~\eqref{muj}) and $v_r:=v_{n,r}^m$ (cf.~\eqref{vr})
 
\begin{minipage}[t]{0.45\textwidth} 
\begin{algorithm}[H]
\caption{Decomposition of $\vec{a}_{3n}$}
{\small
\KwIn{$\vec{a}_{3n}=(a_{3n,1},\ldots, a_{3n,3n})^T$}
\KwOut{$\left\{\begin{array}{ll}\vec{a}_{n}=&(a_{n,1},\ldots, a_{n,n})^T\\
\vec{b}_{n}=&(b_{n,1},\ldots, b_{n,2n})^T
\end{array}\right.$}\vspace{.2cm}
\begin{algorithmic}[1] 
\STATE Compute $\vec{t}=\T_{3n}(\vec{a}_{3n})$
\FOR{$r=0:(n-m)$}\vspace{.1cm}
\STATE $w_{r}=t_r$
\ENDFOR 
\FOR{$r=(n-m+1):(n-1)$}\vspace{.1cm}
\STATE $\displaystyle w_{r}=\mu_r t_r -\mu_{2n-r}t_{2n-r}$
\ENDFOR
\STATE $
\vec{a}_n=\T_n^{\prime}(\vec{w})$ 
\vspace{.2cm}
\FOR{$r=n:n+m-1$}\vspace{.1cm}
\STATE $\displaystyle u_r=\frac{\mu_r t_{2n-r} +\mu_{2n-r}t_r}{\sqrt{v_r}}$ \vspace{.1cm}
\ENDFOR
\FOR{$r=(n+m):(3n-m)$}
\STATE $\displaystyle u_r=\frac{t_r}{\sqrt{v_r}}$\vspace{.2cm}
\ENDFOR
\FOR{$r=(3n-m+1):(3n-1)$}
\STATE $u_r=t_r\sqrt{v_r}$\vspace{.2cm}
\ENDFOR
\STATE $\vec{b}_n=\Sigma_n^{\prime}(\vec{u})$  
\end{algorithmic}
}
\end{algorithm}
\end{minipage}%
\hspace{.2cm}
\begin{minipage}[t]{0.45\textwidth}
\begin{algorithm}[H]
\caption{Reconstruction of $\vec{a}_{3n}$}
{\small
\KwIn{$\left\{\begin{array}{ll}\vec{a}_{n}=&(a_{n,1},\ldots, a_{n,n})^T\\
\vec{b}_{n}=&(b_{n,1},\ldots, b_{n,2n})^T
\end{array}\right.$}
\KwOut{$\vec{a}_{3n}=(a_{3n,1},\ldots, a_{3n,3n})^T$}\vspace{.2cm}
\begin{algorithmic}[1] 
\STATE Compute $\vec{\alpha}=\T_{n}(\vec{a}_{n})$
\STATE Compute $\vec{\beta}=\Sigma_{n}(\vec{b}_{n})$
\FOR{$r=n:3n-1$}\vspace{.1cm}
\STATE $\displaystyle \beta_r=\frac{\beta_r}{\sqrt{v_r}}$
\ENDFOR 
\FOR{$r=0:(n-m)$}\vspace{.1cm}
\STATE $t_r=\alpha_r$
\ENDFOR
\FOR{$r=(n-m+1):(n-1)$}\vspace{.1cm}
\STATE $t_r=\alpha_r \ \mu_r+\beta_{2n-r}\ \mu_{2n-r}$\vspace{.1cm}
\ENDFOR
\STATE $t_n=\beta_n$
\FOR{$r=(n+1):(n+m-1)$}
\STATE $t_r=\beta_r\ \mu_{2n-r} - \alpha_{2n-r}\ \mu_r$\vspace{.2cm}
\ENDFOR
\FOR{$r=(n+m):(3n-m)$}
\STATE $t_r=\beta_r$\vspace{.2cm}
\ENDFOR
\FOR{$r=(3n-m+1):(3n-1)$}
\STATE $t_r=\beta_r\ v_r$\vspace{.2cm}
\ENDFOR
\STATE $\vec{a}_{3n}=\T_{3n}^{\prime}(\vec{t})$  
\end{algorithmic}
}
\end{algorithm}
\end{minipage}
\vspace{1cm}\newline
Regarding the computational cost of the previous algorithms, we remark that all the transforms in \eqref{T}-\eqref{S-1}  can be efficiently computed using $\bigO(n\log n)$ flops by means of the well--known {\em discrete cosine transform} (DCT) and its inverse (IDCT),  in Matlab known  as DCT of type 2 and 3, respectively.  We recall that for any pair of vectors $\vec{v}, \vec{u}\in\RR^N$ they are defined as follows
\begin{itemize}
\item
$\vec{u}=DCT (\vec{v})$ has entries
$\displaystyle
u_r=\sqrt{\frac \pi N}\sum_{k=1}^N 
v_k \ p_{r} (x_k^N),
\quad r=0,\ldots,N-1$,
\item
$\vec{u}=IDCT (\vec{v})$ has entries
$\displaystyle 
u_k=\sqrt{\frac \pi N}\sum_{r=0}^{N-1} v_r p_{r}(x_k^N),
\quad k =1,\ldots,N$.
\end{itemize}
Using these transforms, we have the following algorithms for computing the transforms $\Sigma_n$ and $\Sigma^{\prime}_n$

\begin{minipage}[t]{0.45\textwidth}
\begin{algorithm}[H]
\caption{Compute $\vec{s}=\Sigma_n(\vec{u})$}
{\small
\KwIn{$\vec{u}=(u_1,\ldots, u_{2n})^T$}
\KwOut{$\vec{s}=(s_n,\ldots,s_{3n-1})^T$}\vspace{.2cm}
}
\begin{algorithmic}[1] 
\FOR{$k=1:n$}\vspace{.1cm}
\STATE $w_{3k-2}=u_{2k-1}$
\STATE $w_{3k-1}=0$
\STATE $w_{3k}=u_{2k}$\vspace{.1cm}
\ENDFOR \vspace{.2cm}
\STATE $
\vec{x}= DCT(\vec{w})$
\vspace{.2cm}
\STATE $s_n=x_n$
\FOR{$k=n+1:2n-1$}\vspace{.1cm}
\STATE $s_k=\frac{x_k+x_{2n-k}}{\sqrt{2}}$\vspace{.1cm}
\ENDFOR\vspace{.1cm}
\STATE $s_{2n}=\frac{x_{2n}+\sqrt{2}x_0}{\sqrt{3}}$
\FOR{$k=2n+1:3n-1$}\vspace{.2cm}
\STATE $s_k=\sqrt{\frac 32} x_k$\vspace{.2cm}
\ENDFOR
\end{algorithmic}
\end{algorithm}
\end{minipage}%
\hspace{.2cm}
\begin{minipage}[t]{0.45\textwidth}
\begin{algorithm}[H]
{\small
\caption{Compute $\vec{u}=\Sigma_n^{\prime}(\vec{s})$}
\KwIn{$\vec{s}=(s_n,\ldots, s_{3n-1})^T$}
\KwOut{$\vec{u}=(u_1,\ldots,u_{2n})^T$}\vspace{.2cm}
\begin{algorithmic}[1]
{
\STATE $w_0=\sqrt{\frac 23}s_{2n}$
\FOR{$k=1:n-1$}
\STATE $w_k=s_{2n-k}/\sqrt{2}$
\ENDFOR
\STATE $w_n=s_n$
\FOR{$k=n+1:2n-1$} 
\STATE $w_k=s_k/\sqrt{2}$
\ENDFOR\vspace{.1cm}
\STATE $w_{2n}=\frac{s_{2n}}{\sqrt{3}}$\vspace{.2cm}
\FOR{$k=2n+1:3n-1$}
\STATE $w_k=\sqrt{\frac 32} s_k$
\ENDFOR \vspace{.2cm}
\STATE $\displaystyle
\vec{x}= IDCT(\vec{w})$\vspace{.1cm}
\FOR{$k=1:n$}
\STATE $u_{2k-1}=x_{3k-2}$
\STATE $u_{2k}=x_{3k}$
\ENDFOR
}
\end{algorithmic}
}
\end{algorithm}
\end{minipage}
\vspace{.5cm}\newline
Moreover, using the parameters $\nu_r:=\nu_{n,r}^n$ defined in \eqref{nur},
the transforms $\T_n$ and $\T^\prime_n$ can be computed as follows 

\begin{minipage}[t]{0.45\textwidth}
\begin{algorithm}[H]
\caption{Compute $\vec{t}=\T_n(\vec{u})$}
{\small
\KwIn{$\vec{u}=(u_1,\ldots, u_{2n})^T$}
\KwOut{$\vec{t}=(t_0,\ldots,t_{n-1})^T$}\vspace{.2cm}
}
\begin{algorithmic}[1] 
\STATE $\vec{t}= DCT(\vec{u})$
\FOR{$r=n-m+1:n-1$}\vspace{.1cm}
\STATE $t_r=t_r/\sqrt{\nu_r}$\vspace{.1cm}
\ENDFOR
\end{algorithmic}
\end{algorithm}
\end{minipage}
\hspace{.2cm}
\begin{minipage}[t]{0.45\textwidth}
\begin{algorithm}[H]
{\small
\caption{Compute $\vec{u}=\T_n^{\prime}(\vec{t})$}
\KwIn{$\vec{t}=(t_0,\ldots, t_{n-1})^T$}
\KwOut{$\vec{u}=(u_1,\ldots,u_{n})^T$}\vspace{.2cm}
\begin{algorithmic}[1]
\FOR{$r=0:n-1$}\vspace{.1cm}
\STATE $t_r=t_r/\sqrt{\nu_r}$
\vspace{.1cm}
\ENDFOR
\STATE $\displaystyle
\vec{u}= IDCT(\vec{t})$\vspace{.1cm}
\end{algorithmic}
}
\end{algorithm}
\end{minipage}

\begin{figure}[!thb]%
\begin{center}
\includegraphics[scale = 0.28]{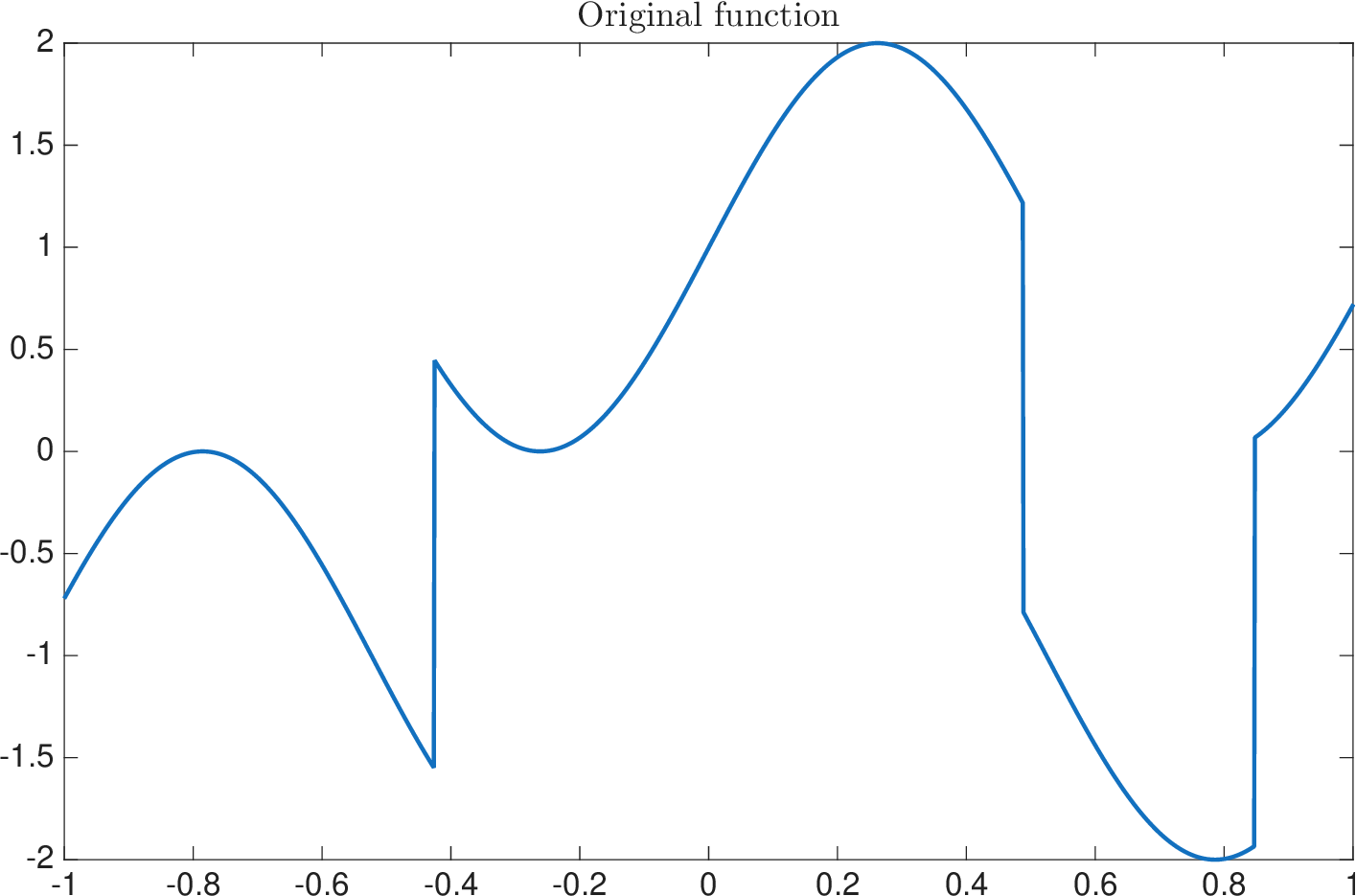}
\hspace{0.2cm}
\includegraphics[scale = 0.28]{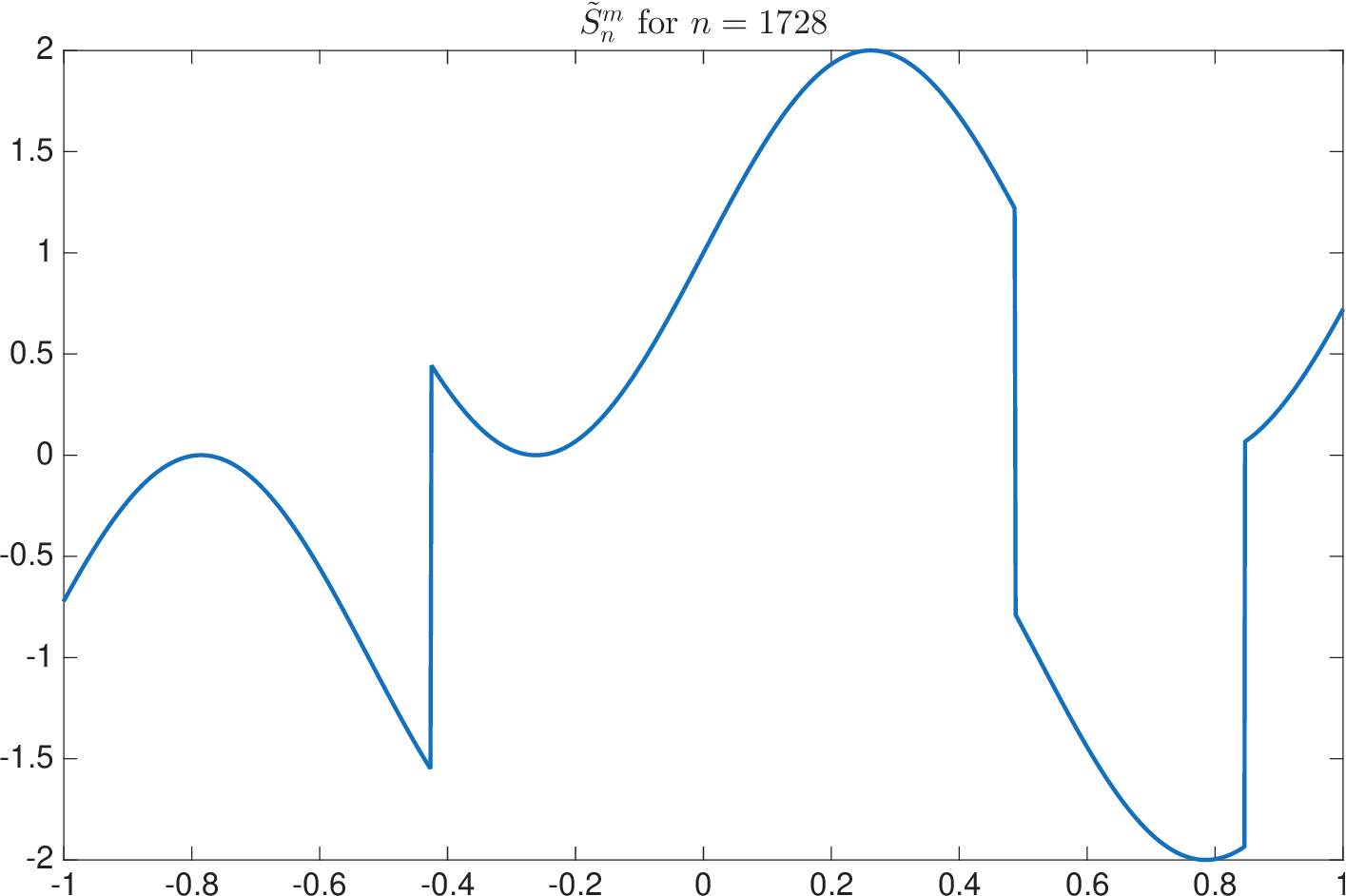}\\[0.3cm]
\includegraphics[scale = 0.28]{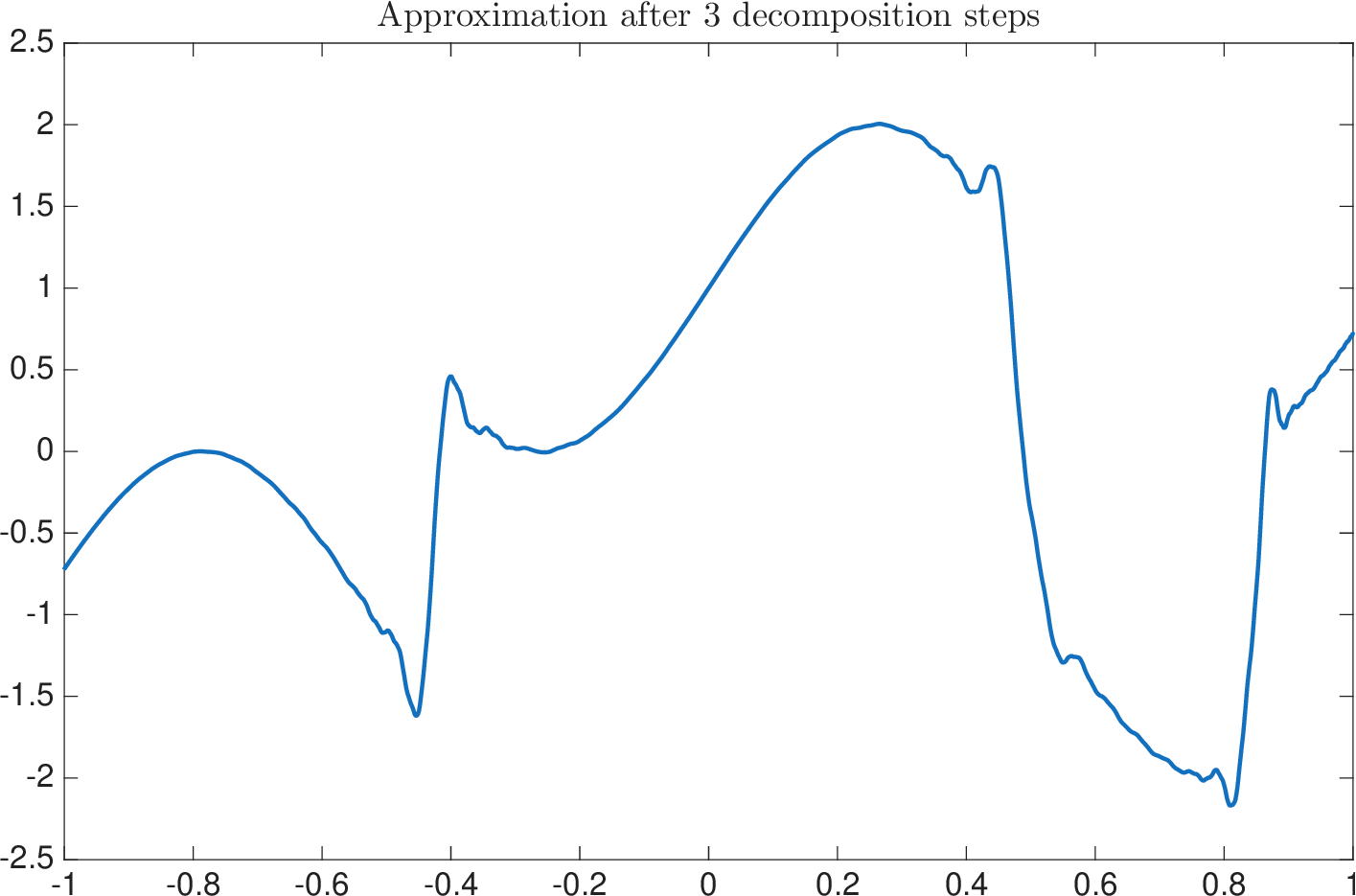}
\hspace{0.2cm}
\includegraphics[scale = 0.28]{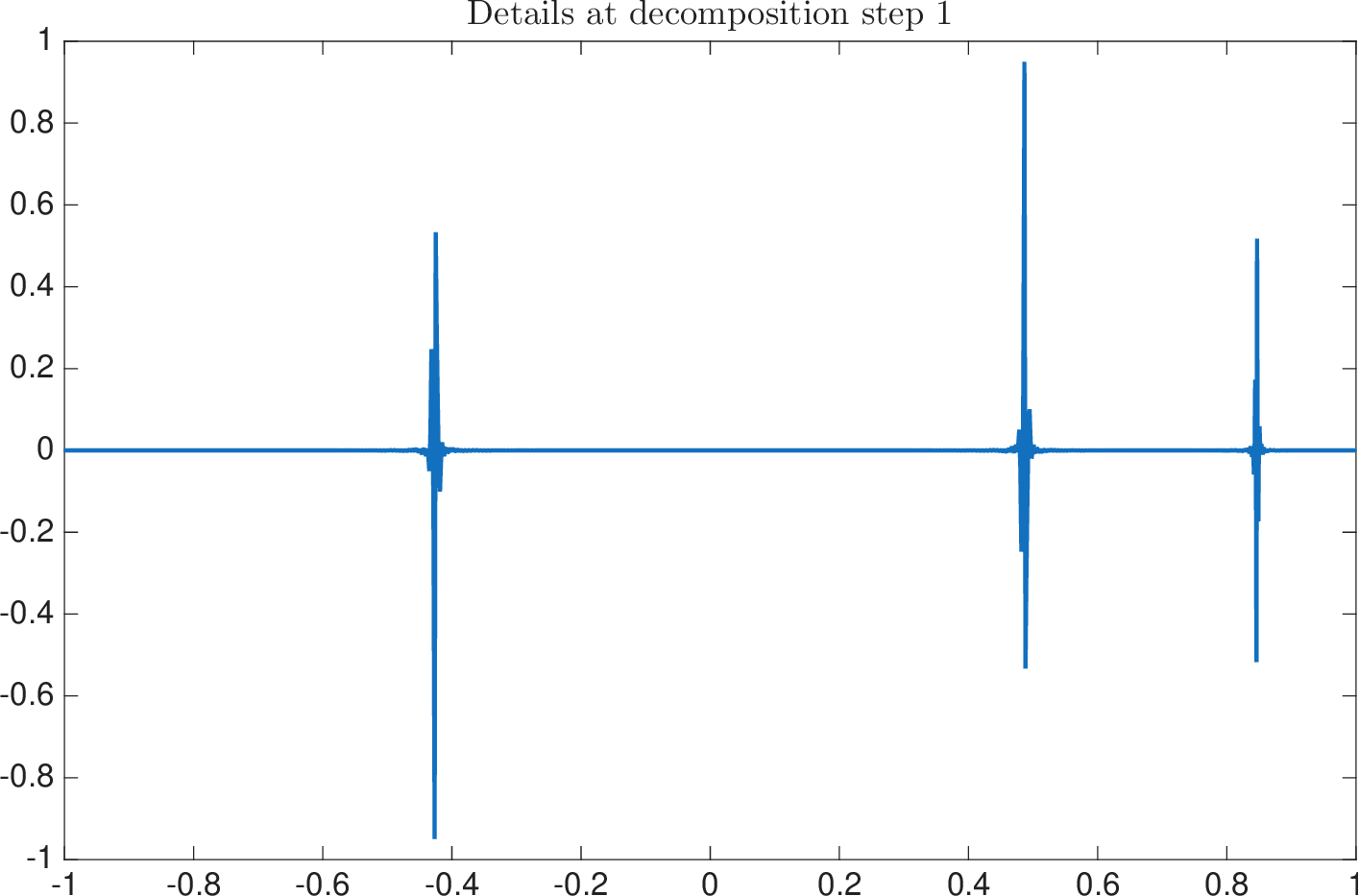}\\[0.3cm]
\includegraphics[scale = 0.28]{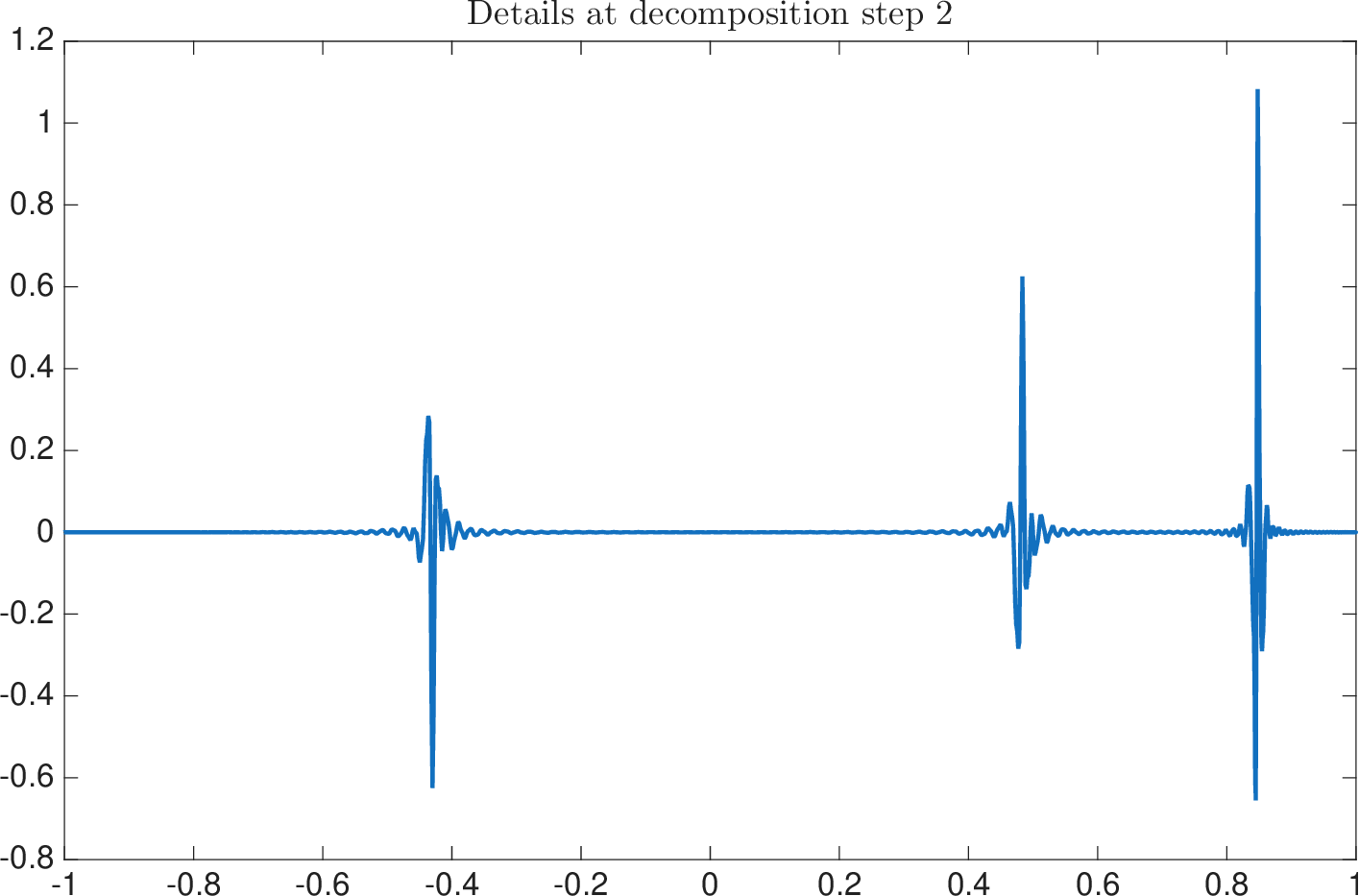}
\hspace{0.2cm}
\includegraphics[scale = 0.28]{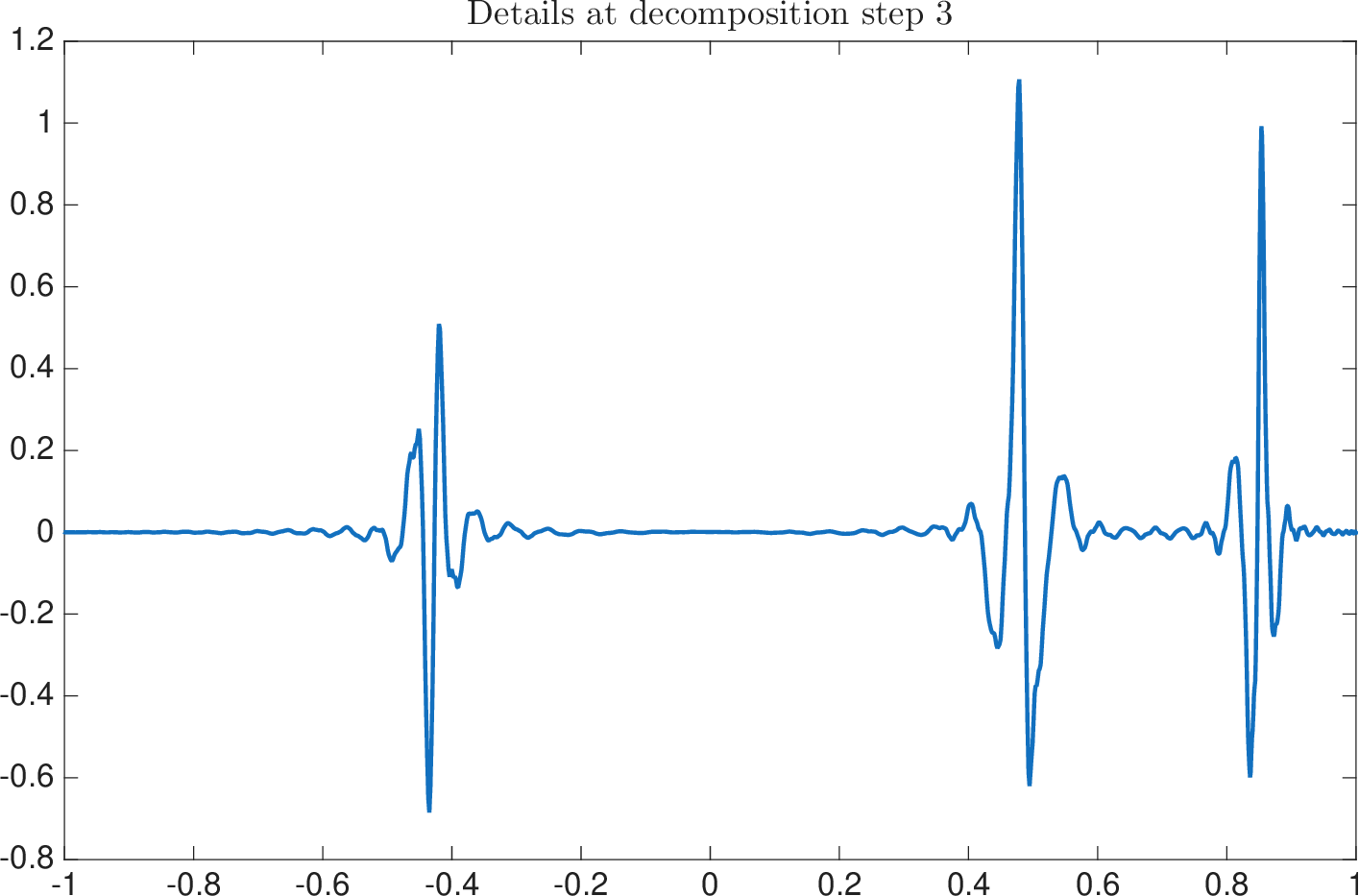}
\end{center}
\caption{ An example of three decomposition steps applied to the approximation $\tilde{S}_n^m f(x)$ of the function
$f(x) = \sin(6x)+\textrm{sign}(\sin(x+\exp(2x)))$ with $n=1728$ and $m =\lfloor 0.7 n\rfloor$.
\label{fig_wavelet_decomposition} }
\end{figure}

In conclusion, we present an example of three decomposition steps for the function
\[
f(x) = \sin(6x) + \mathrm{sign}(\sin(x + \exp(2x))),
\]
using \( m = \lfloor 0.7\, n \rfloor \) at each level  \( n \).
Starting from \( n = 64 \cdot 3^3 = 1728 \), we first project \( f \) onto \( \V_n^m \) by computing the polynomial approximation \( f_{1728} := \tilde{S}_n^m f \) based on the samples of $f$ on \( X_n \).
Then we iteratively decompose \( f_{1728} \) as follows
\[
\begin{aligned}
f_{1728} &= f_{576} + g_{1152} \\
         &= f_{192} + g_{384} + g_{1152} \\
         &= f_{64} + g_{128} + g_{384} + g_{1152}.
\end{aligned}
\]
In Figure \ref{fig_wavelet_decomposition} the function $f$, its initial approximation $f_{1728}$, the lower degree approximation $f_{64}$, and the details $g_{128}$, $ g_{384}$ and $g_{1152}$ are plotted.

\section{Comparison with the interpolating case}\label{sec-experiments}

In Section \ref{sec-approx-oper} we have seen that the Fourier projection $S_n^mf$ associated with the orthonormal VP scaling functions \eqref{sca-ort}, as well as its discrete counterpart $\tilde S_n^mf$, gives a well-conditioned and near best approximation of $f$ in the sampling space $\V_n^m$.

\begin{figure}[!htb]%
\begin{center}
\includegraphics[scale = 0.50]{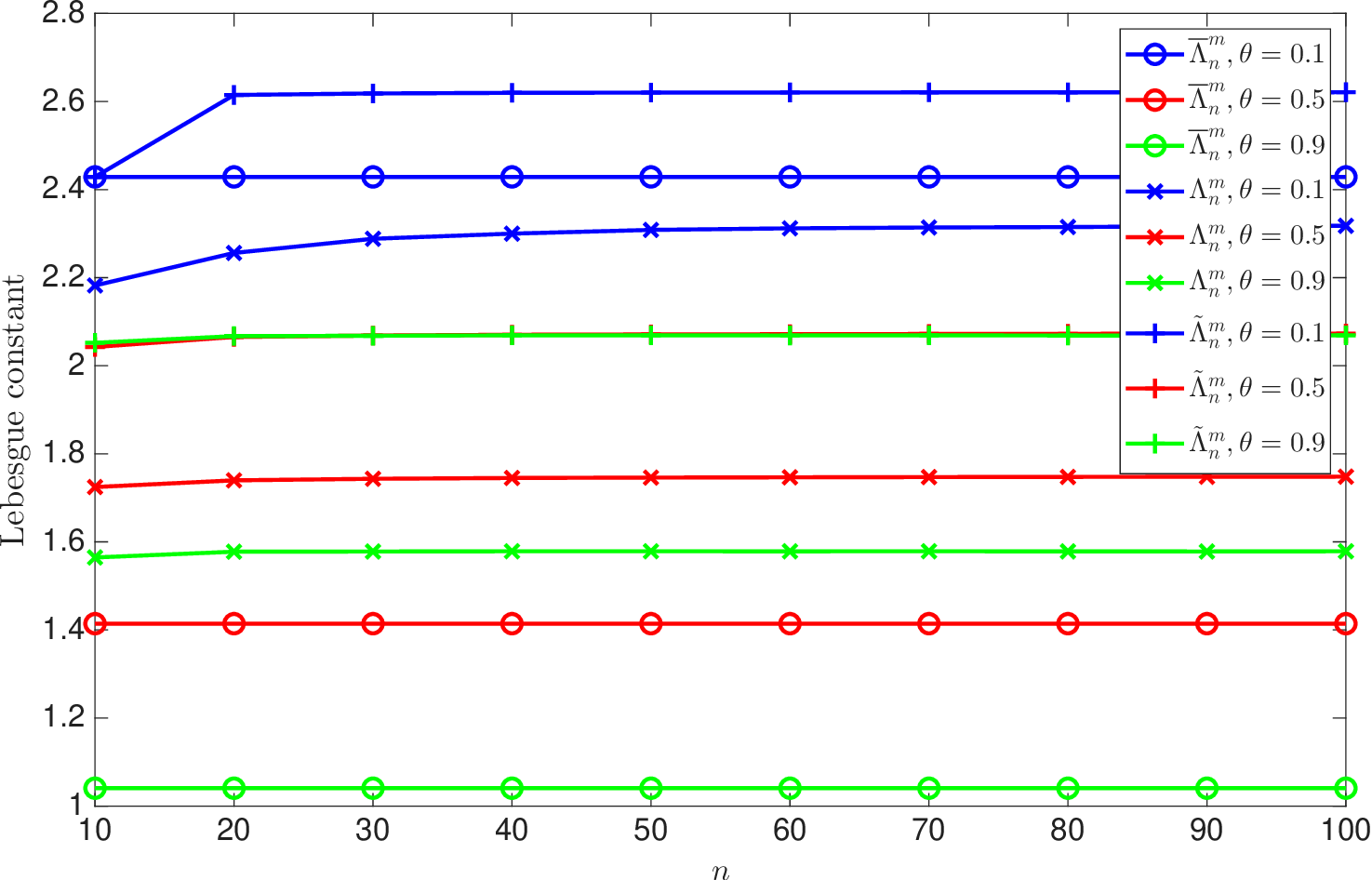}
\vspace{-.5cm}
\end{center}
\caption{ Lebesgue constants
$\overline{\Lambda}_n^m$ (\textendash$\circ$\textendash), $\Lambda_n^m$ (\textendash{\footnotesize $\times$}\textendash), and
$\tilde\Lambda_m^n$ (\textendash{\footnotesize $|$}\textendash), for $n = 10,20,\ldots,100$ and $m=\lfloor \theta n\rfloor$, with $\theta = 0.1 \mbox{ (blue)}, 0.5 \mbox{ (red)}, 0.9 \mbox{ (green)}$.
\label{fig_LC2} }
\end{figure}

We recall that similar properties have already been proved in the literature  by another projection on $\V_n^m$ based  on the interpolating VP scaling functions \eqref{sca}, and defined by 
\begin{equation}\label{VP}
V_n^mf(x):=\sum_{i=1}^n f(x_i^n)\Phi_{n,i}^m(x), \qquad \forall x\in [-1,1].
\end{equation}
This polynomial is known in the literature as the VP interpolation polynomial  of $f$ and its approximation properties  have been investigated in detail in several papers (see, e.g.,  \cite{OT-APNUM21, OT-DRNA21} and the references therein). Here, we only recall that the results stated in Theorem \ref{th-tildeSn}  continue to hold by replacing  $\tilde S_n^mf$  with the polynomial $V_n^mf$, which, in particular, satisfies
\begin{equation}\label{LC-VP}
\|V_n^mf\|_\infty\le \overline{\Lambda}_n^m \max_{1\le k\le n} |f(x_k^n)|, \qquad \overline{\Lambda}_n^m:= \max_{|x|\le 1}\sum_{k=1}^n|\Phi_{n,k}^m(x)|.
\end{equation}
Like the Lebesgue constants $\Lambda_n^m$ in \eqref{def-LC} and $\tilde\Lambda_n^m$ in \eqref{def-LC1}, also the Lebesgue constants $\overline{\Lambda}_n^m$ in \eqref{LC-VP} are uniformly bounded w.r.t.\ $n$ provided that $m=\lfloor \theta n\rfloor$ with $\theta\in ]0,1[$ arbitrarily fixed (see, e.g., \cite[Thm 3.1(e)]{TB-wave}. In Figure \ref{fig_LC2} their behavior is compared for three different values of $\theta$.

We note that, by \eqref{sca-ort1-inv}, the VP interpolation polynomial \eqref{VP} can be expanded in the new orthonormal VP scaling basis  as follows  
\begin{equation}\label{VP-1}
V_n^mf(x)=    \sum_{k=1}^n\left[\left(\frac\pi n\right )^\frac 32\sum_{h=1}^n f(x_h^n) \sum_{r=0}^{n-1}\sqrt{\nu_{n,r}^m}p_r(x_k^n)p_r(x_h^n)\right]\tilde\varphi_{n,k}^m(x).
\end{equation}
Comparing \eqref{Sn-tilde2} and \eqref{VP-1} we see the difference between the near best approximation polynomials $\tilde S_n^mf$ and  $V_n^mf$, both based on the same function values. However, by \eqref{sca-inter}, the polynomial $V_n^mf$ also satisfies  the interpolation property
\begin{equation}\label{VP-inter}
 V_n^mf(x_k^n)=f(x_k^n), \qquad k=1,\ldots,n , \qquad 0<m<n  .
\end{equation}

\begin{figure}[!b]%
\begin{center}
\includegraphics[scale = 0.40]{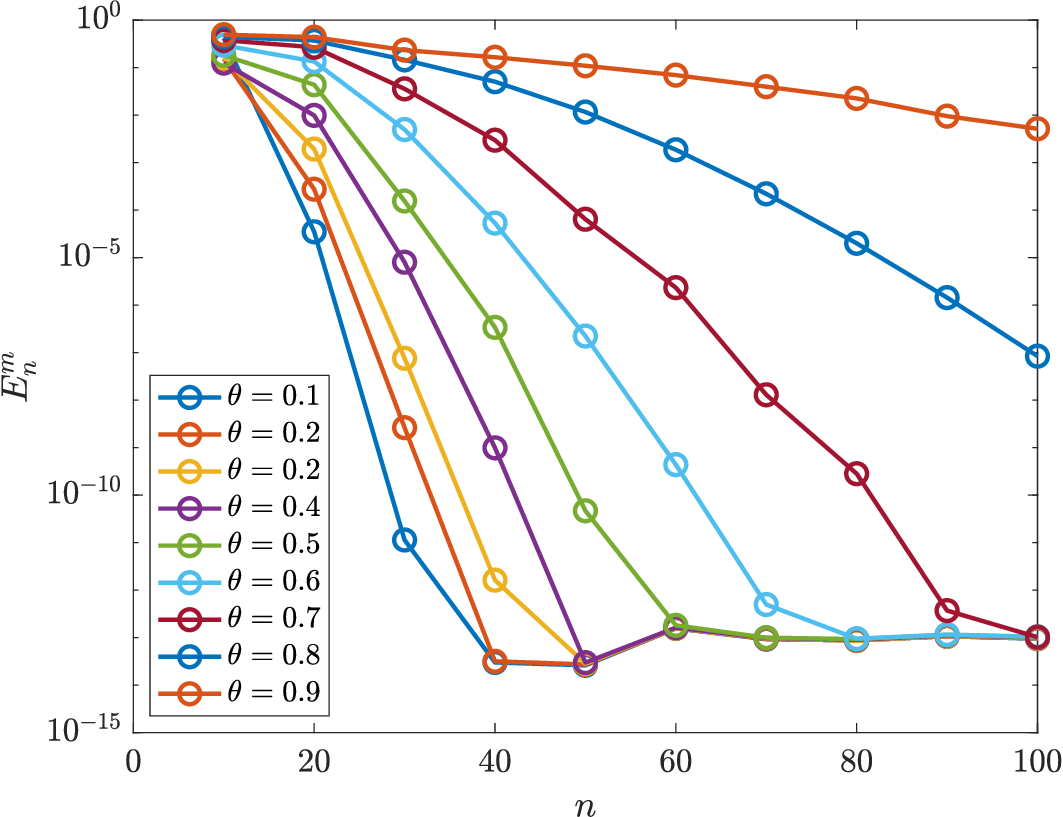}
\hspace{0.2cm}
\includegraphics[scale = 0.40]{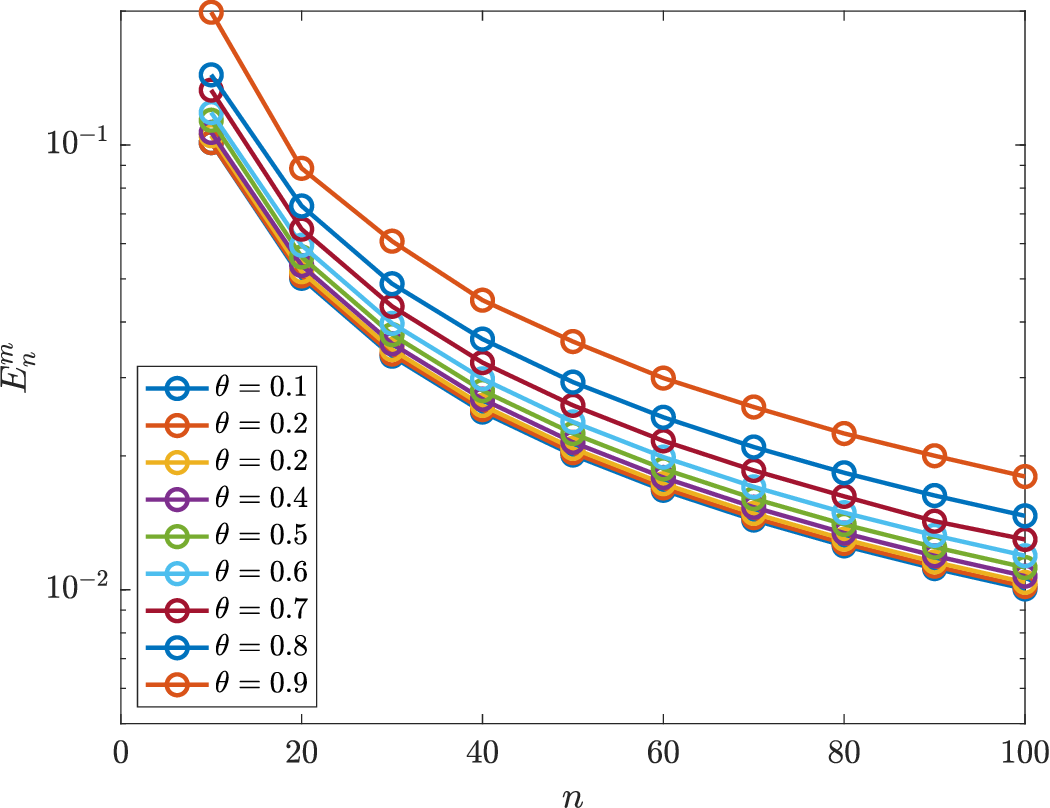}
\\
\includegraphics[scale = 0.40]{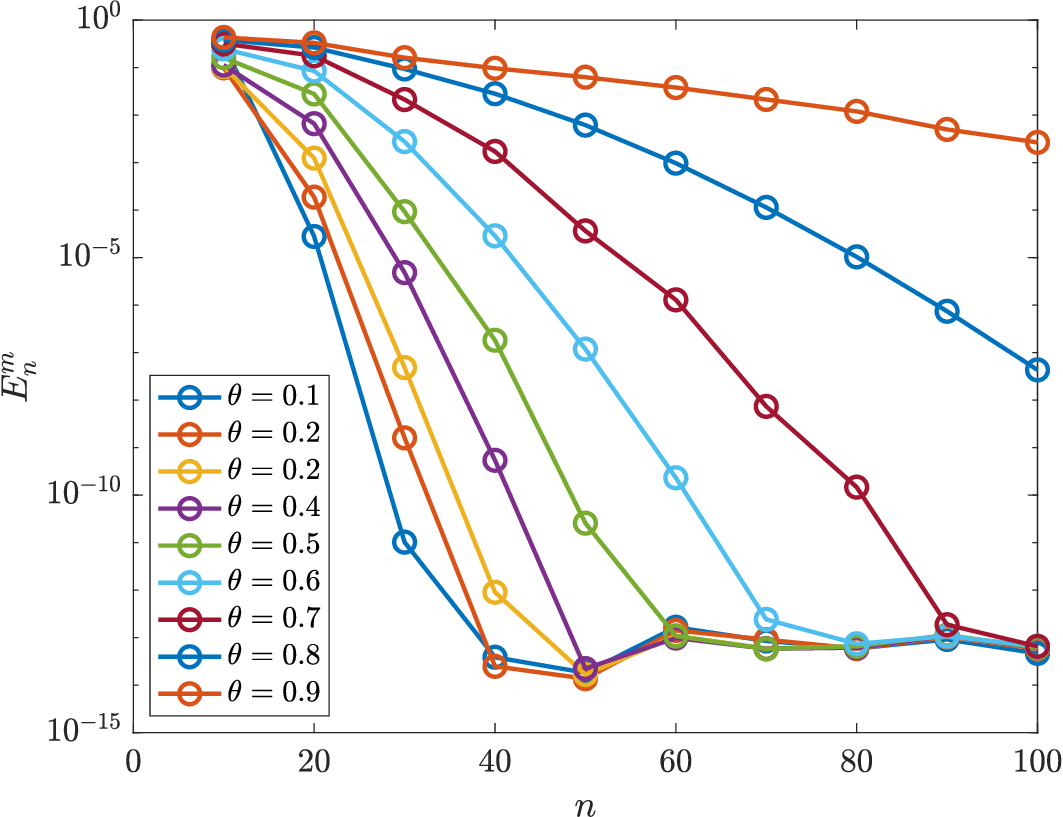}
\hspace{0.2cm}
\includegraphics[scale = 0.40]{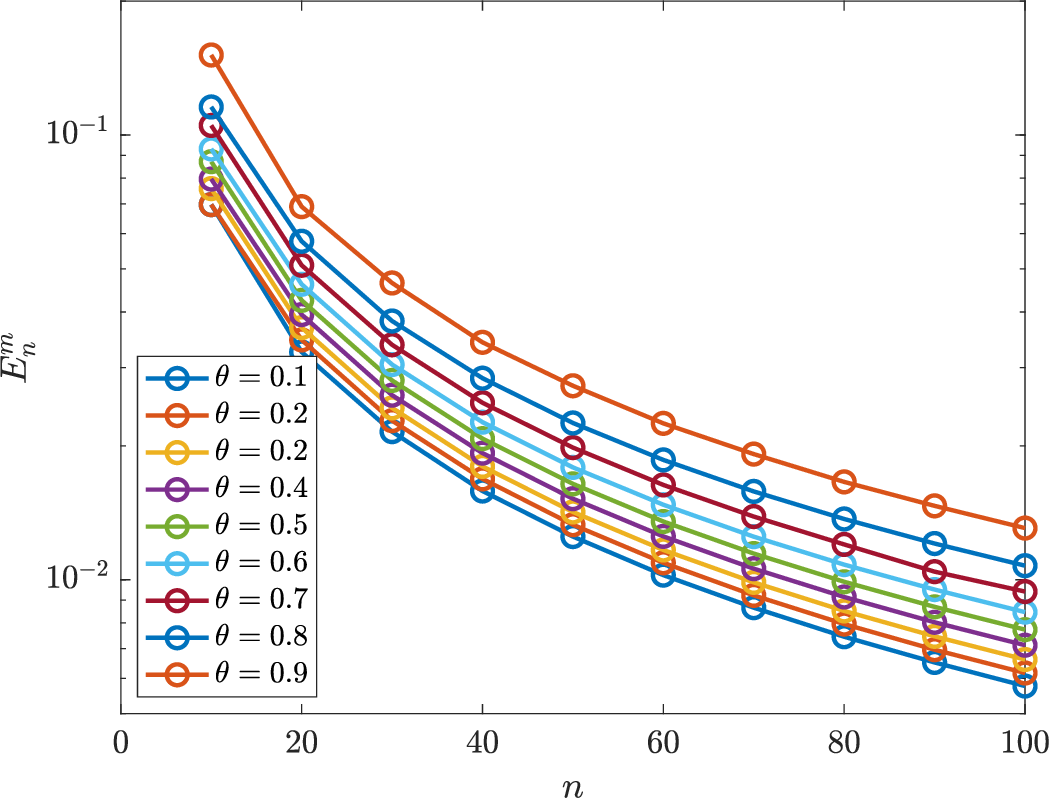}
\\
\includegraphics[scale = 0.40]{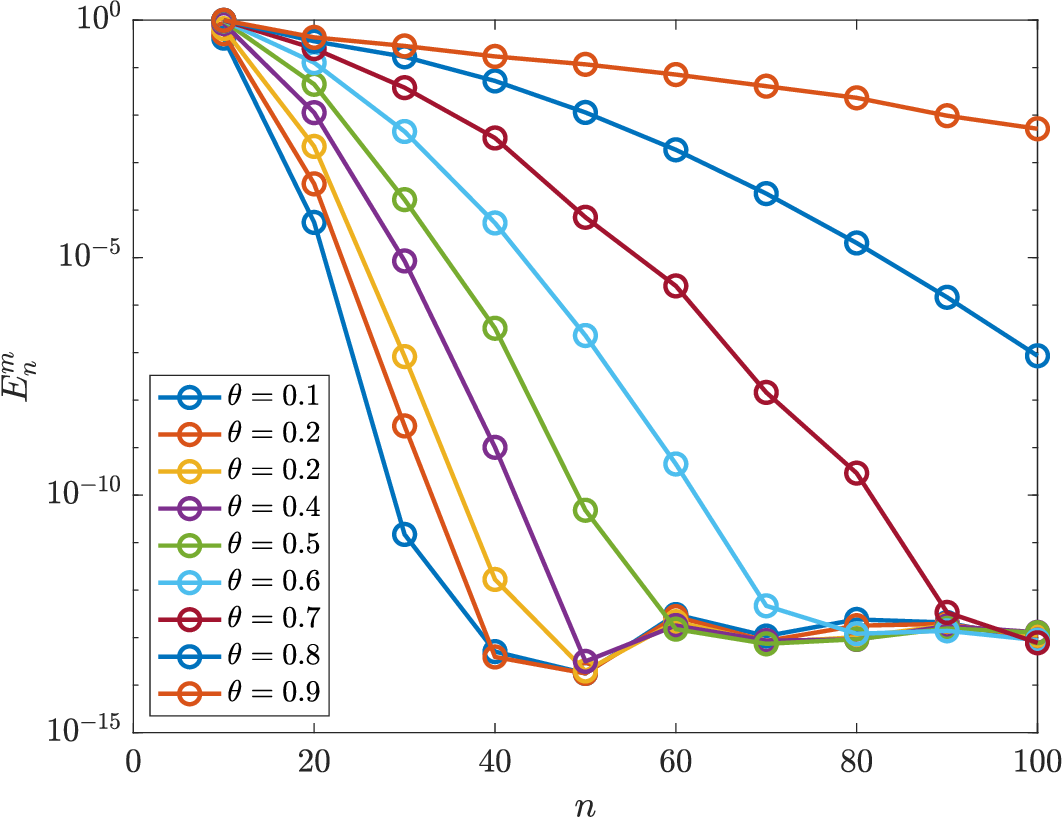}
\hspace{0.2cm}
\includegraphics[scale = 0.40]{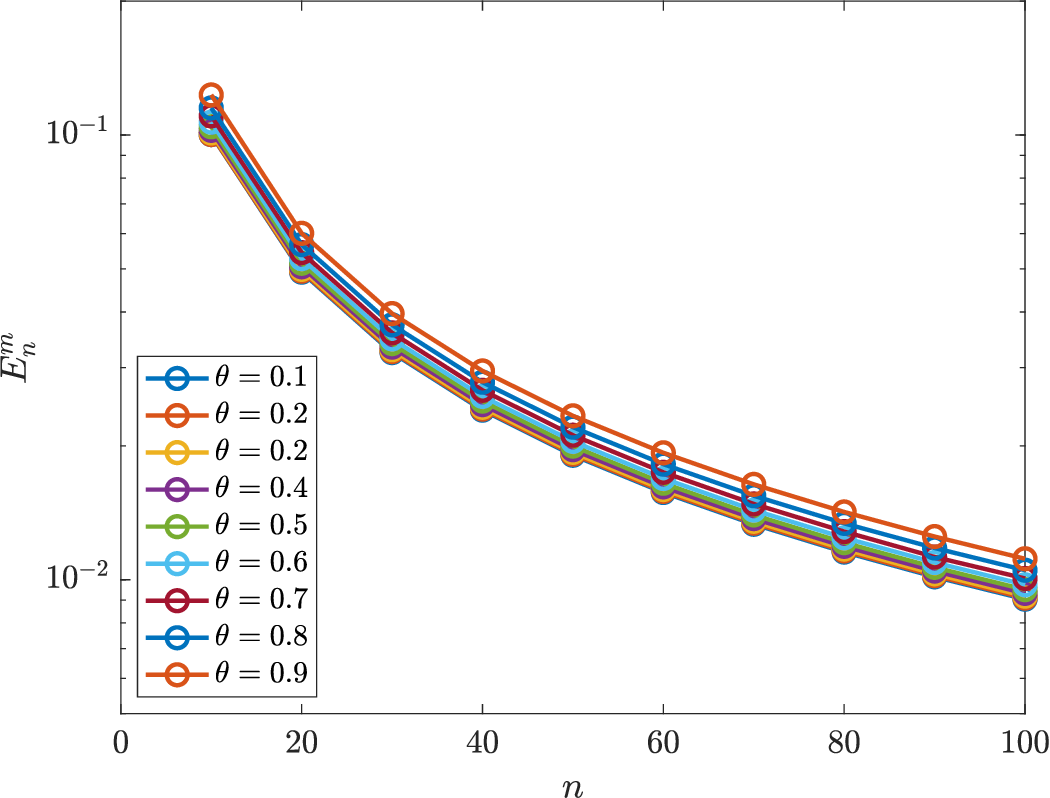}
\end{center}
\caption{ The maximum absolute error $E_n^m$ for the three approximants
$V_n^m f(x)$ (top), $S_n^m f(x)$ (middle) and $\tilde{S}_n^mf(x)$ (bottom) applied to the functions
$f(x)=\sin(x)$ (left) and $f(x)=\mbox{abs}(x)$ (right).\label{fig_E} }
\end{figure}

\begin{figure}[!b]%
\begin{center}
\includegraphics[scale = 0.28]{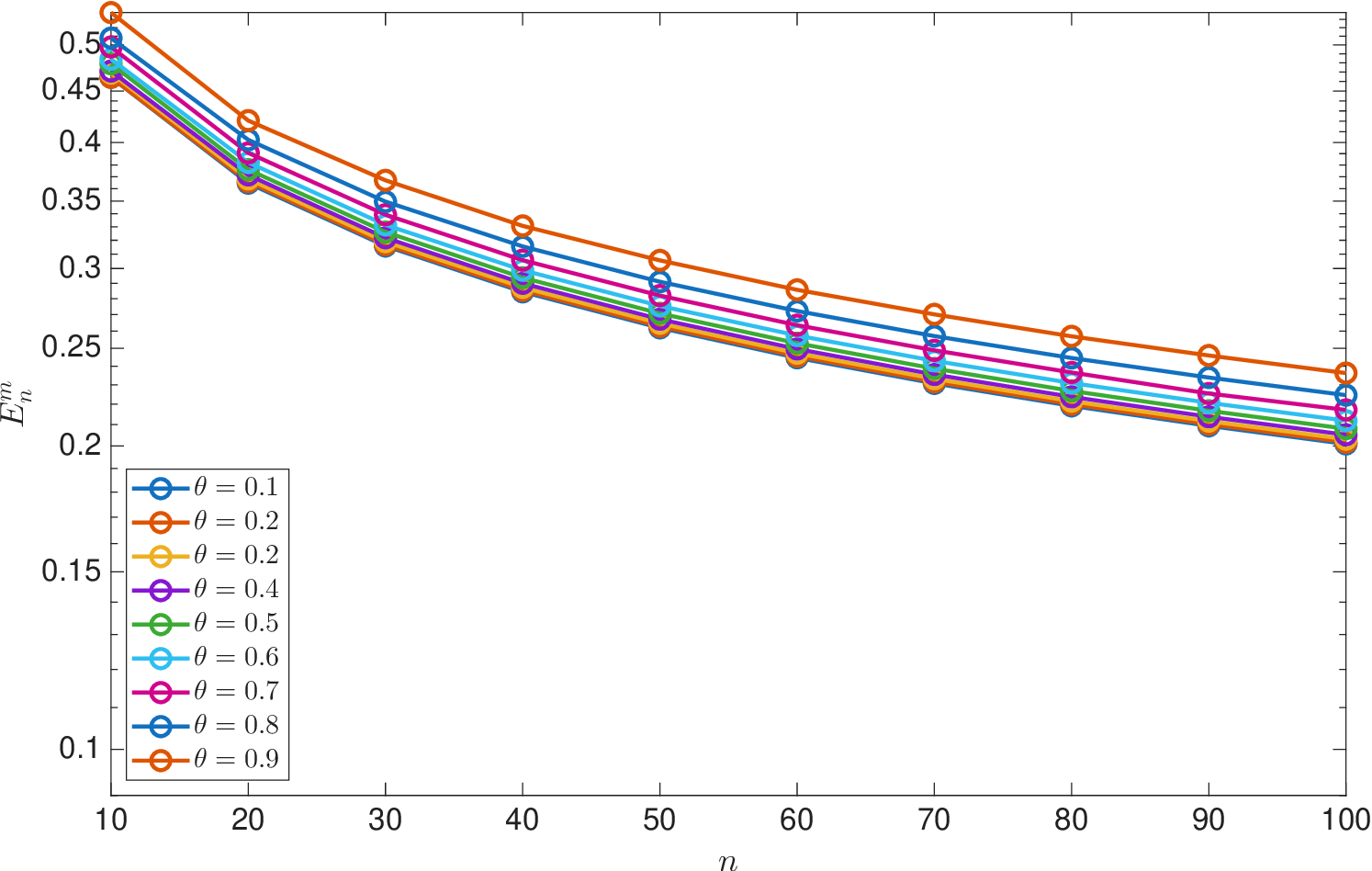}
\hspace{0.2cm}
\includegraphics[scale = 0.28]{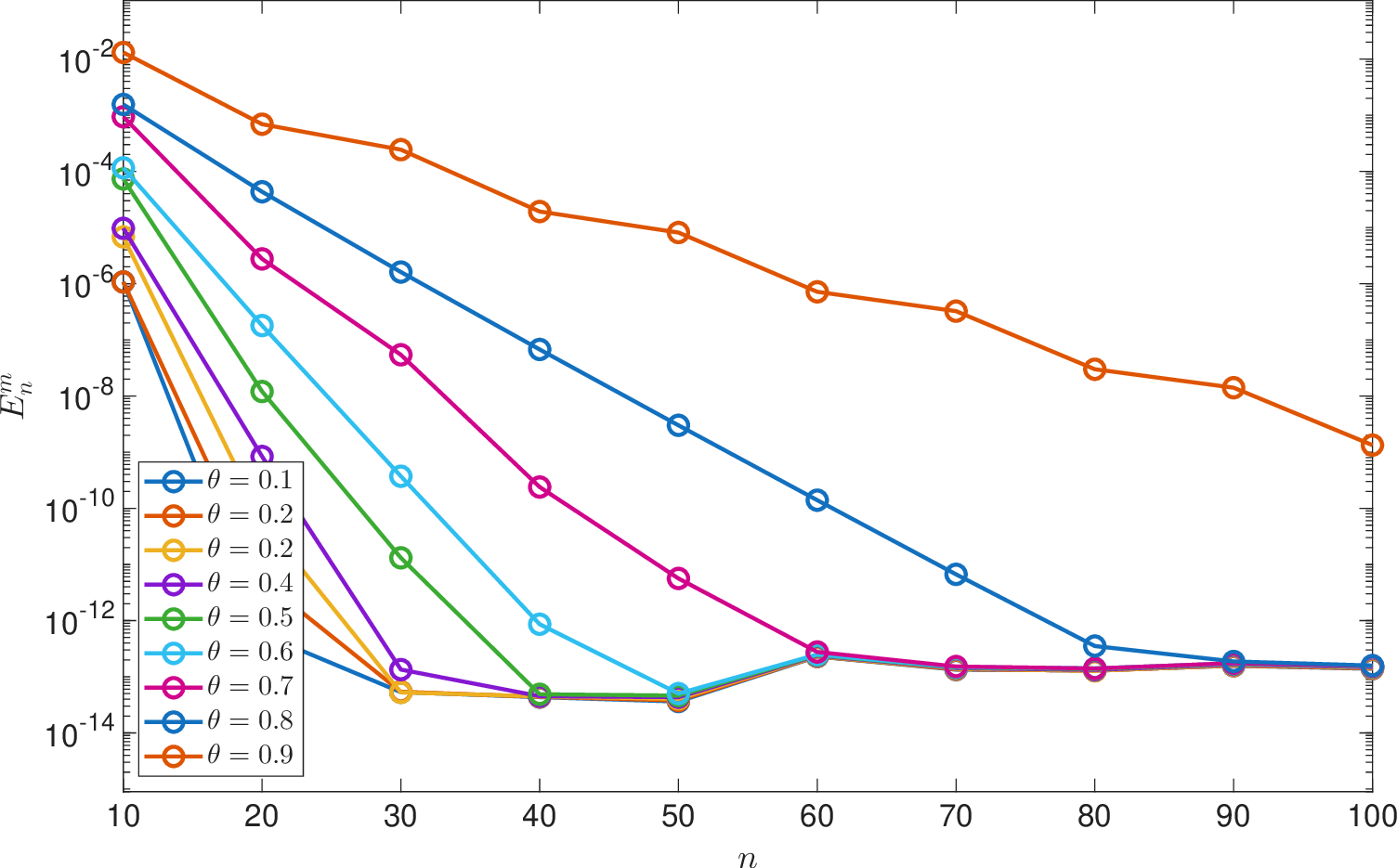}
\\
\includegraphics[scale = 0.28]{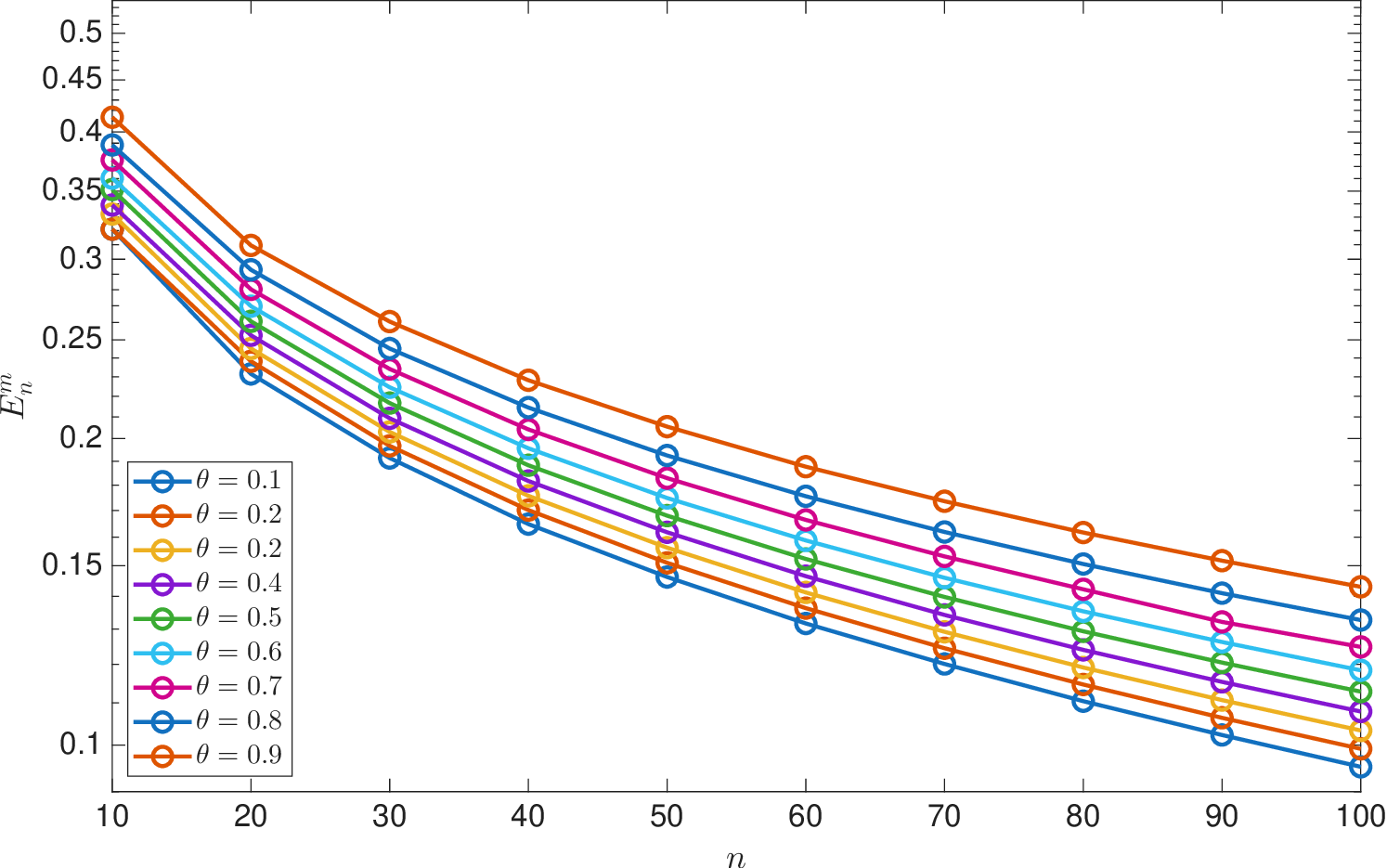}
\hspace{0.2cm}
\includegraphics[scale = 0.28]{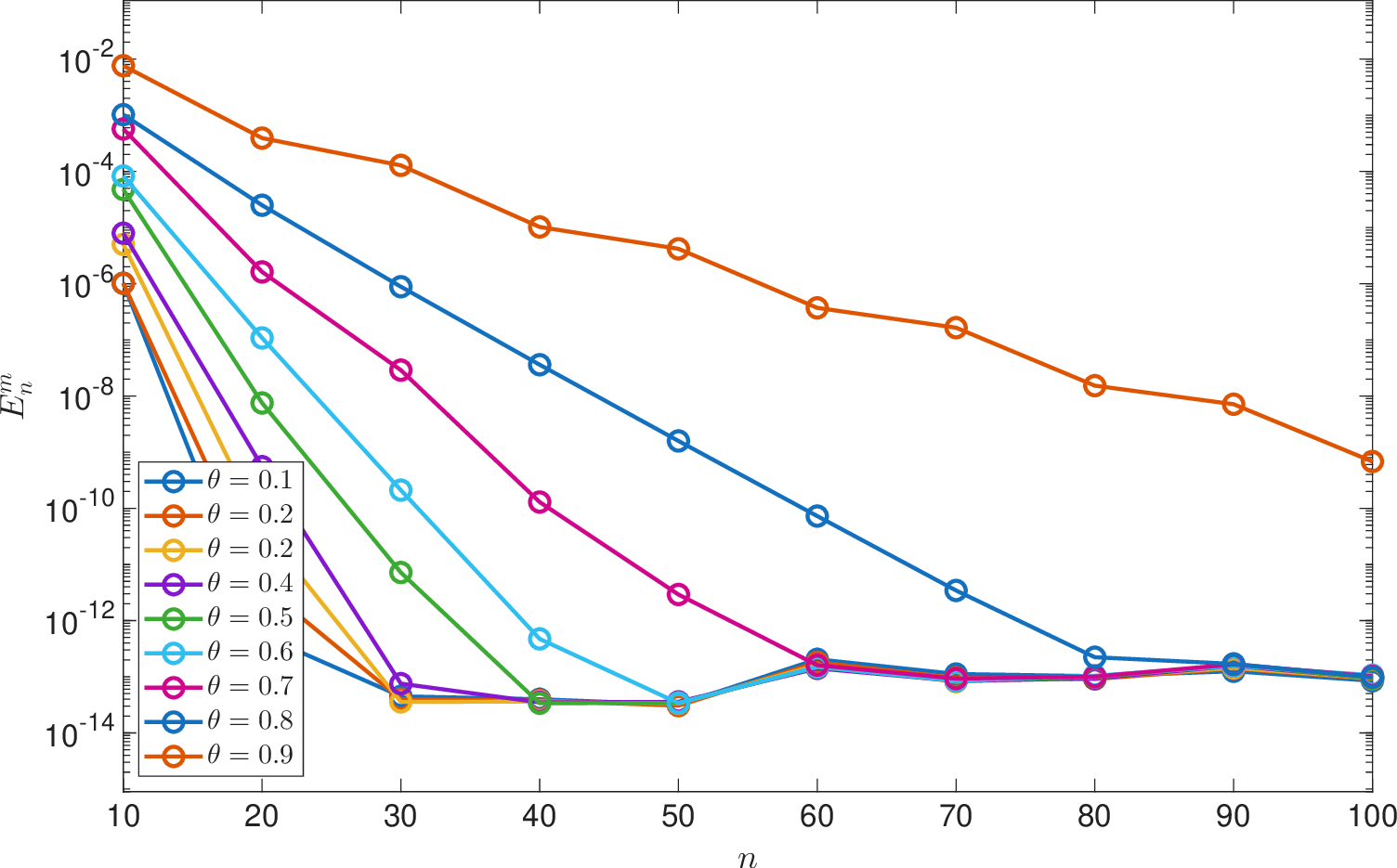}
\\
\includegraphics[scale = 0.28]{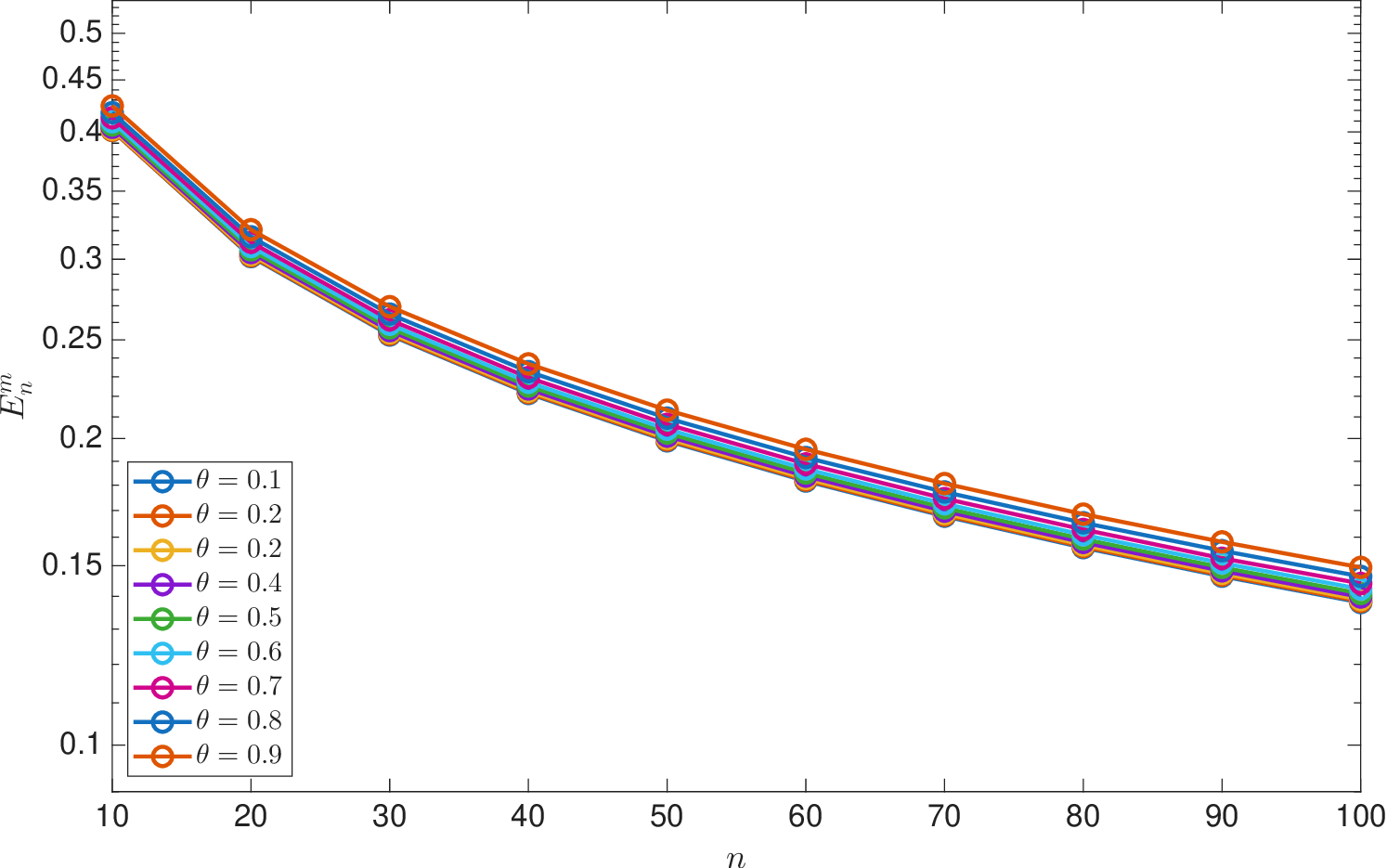}
\hspace{0.2cm}
\includegraphics[scale = 0.28]{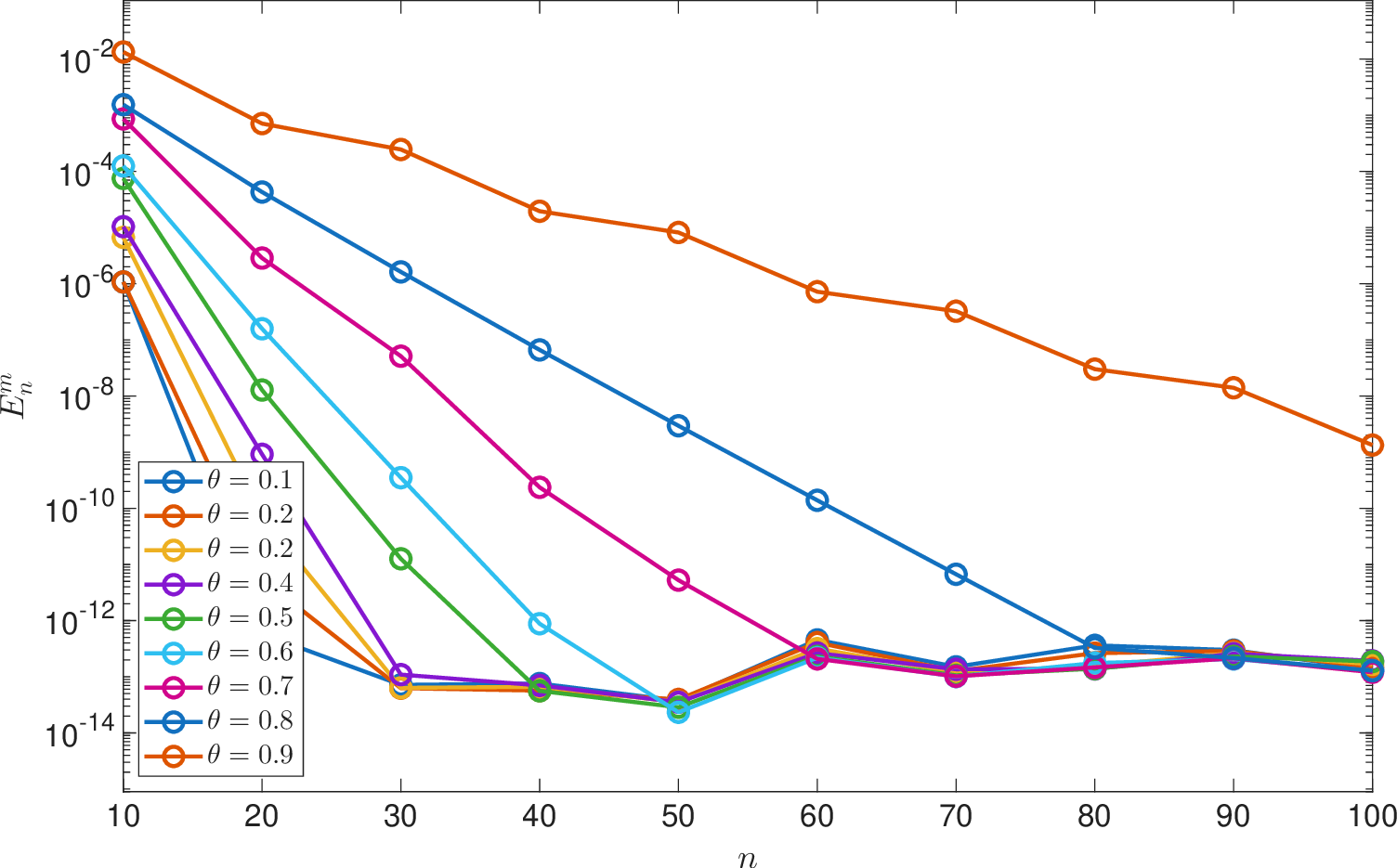}
\end{center}
\caption{ The maximum absolute error $E_n^m$ for the three approximants
$V_n^m f(x)$ (top), $S_n^m f(x)$ (middle) and $\tilde{S}_n^mf(x)$ (bottom) applied to the functions
$f(x)=\mbox{abs}(x)^{0.3}$ (left) and $f(x)=1/(1+0.25x^2)$ (right).\label{fig_EE} }
\end{figure}

To illustrate the approximation properties of the three approximants,
$V_n^m f(x)$, $S_n^m f(x)$ and $\tilde{S}_n^mf(x)$, the approximation
error $E_n^m$ is 
plotted for the functions $\sin(x)$ and $\mbox{abs}(x)$ in Figure \ref{fig_E},
and for the functions $\mbox{abs}(x)^{0.3}$ and $1/(1+0.25x^2)$ in Figure~\ref{fig_EE}.
The approximation error $E_n^m$ is defined as
\begin{equation*}
    E_n^m = \| f_{\mbox{approx}} - f \|_\infty,
\end{equation*}
and we take $m=\lfloor \theta n\rfloor$, with $\theta=0.1,0.2,\ldots, 0.9$ and $n=10,20,\ldots, 100$.

The numerical results show that, for very smooth functions, the errors are very similar, whereas for less regular functions the differences become more pronounced. This is illustrated in the second column of Fig.~\ref{fig_E} and in the first column of Fig.~\ref{fig_EE}, where, for higher degrees, the Fourier projection $S_n^mf$ yields the smallest error, followed by $\tilde S_n^mf$,  whose error is much less sensitive to $\theta$, while $V_n^mf$ gives the largest error.

\section{Proofs}\label{sec-proofs}
\subsection{Proof of Proposition \ref{prop-sca}} 
Since $\tilde\varphi_{n,k}^m$ are defined as a linear combination of the basis polynomials \eqref{q-basis}, it is obvious that $\tilde\varphi_{n,k}^m\in\V_n^m$ for all $k=1,2,\ldots,n$. Hence we only have to prove \eqref{prod-sca}. To this aim we recall that the Darboux kernel
\[
K_n(x,y)=\sum_{s=0}^{n-1} p_s(x)p_s(y), \quad x,y\in [-1,1], \qquad n=0,1,\ldots,
\]
satisfies
\[
\frac \pi n K_n(x_k^n, \ x_h^n)=\delta_{h,k}, \qquad h,k=1,\ldots, n.
\]
Consequently, recalling \eqref{sca-ort}, \eqref{tau}, and \eqref{q-prod}, we get
\begin{eqnarray*}
<\tilde\varphi_{n,k}^m,\ \tilde\varphi_{n,h}^m>_{L^2_w} &=&
\frac \pi n\sum_{r=0}^{n-1}\sum_{s=0}^{n-1}\frac{p_r(x_k^n)}{\sqrt{\nu_{n,r}^m}} \frac{p_s(x_h^n)}{\sqrt{\nu_{n,s}^m}}\
<q_{n,r}^m,\ q_{n,s}^m>_{L^2_w}\\
&=& \frac \pi n \sum_{r=0}^{n-1} p_r(x_k^n) p_r(x_h^n)=\delta_{k,h}. 
\end{eqnarray*}\Proofend

\subsection{ Proof of Proposition \ref{prop-sca1}} 
Firstly, note that by the interpolating property of the scaling functions $\Phi_{n,h}^m$, we have that
\[
f(x)=\sum_{h=1}^nf(x_h^n)\Phi_{n,h}^m(x),\qquad \forall f\in\V_n^m.
\]
On the other hand, we note that
\begin{equation}\label{p-xh}
p_{2n-r}(x_h^n)=-p_r(x_h^n), \qquad h=1,\ldots,n, \qquad r=1,\ldots,n-1,
\end{equation}
and 
\[
\mu_{n,r}^m+\mu_{n,2n-r}^m=1, \qquad n-m<r<n,
\]
imply that
\begin{equation}\label{q-xh}
q_{n,r}^m(x_h^n)=p_r(x_h^n),\qquad h=1,\ldots,n,\qquad r=0,\ldots,n-1.
\end{equation}
Consequently, we get
\begin{equation}\label{q-fi}
q_{n,r}^m(x)=\sum_{h=1}^np_r(x_h^n)\Phi_{n,h}^m(x), \qquad r=0,\ldots,n-1.    
\end{equation}
Using this identity and Definition \ref{def-sca},  we conclude that
\[
\tilde\varphi_{n,k}^m(x)=
\sqrt{\frac\pi n}\ \sum_{r=0}^{n-1}\frac{p_r(x_k^n)}{\sqrt{\nu_{n,r}^m}}\  q_{n,r}^m(x)= \sqrt{\frac\pi n}\sum_{h=1}^n\left[\sum_{r=0}^{n-1}\frac{p_r(x_k^n)}{\sqrt{\nu_{n,r}^m}}\  p_{r}(x_h^n)\right]\Phi_{n,h}^m(x)
\]
holds for $k=1,\ldots,n$.\Proofend

\subsection{Proof of Theorem \ref{th-LC}}
Firstly, we note that by \cite[Thm. 3.1]{TB-wave} we have
\begin{equation}\label{LC-int}
 \sup_{|x|\le 1}\sum_{k=1}^n|\Phi_{n,k}^m(x)|\le \C   
\end{equation}
with $\C>0$ independent of $n$.

Consequently, by \eqref{def-LC}, \eqref{sn-q}, and \eqref{q-fi}, we get
\begin{eqnarray*}
\Lambda_n^m &=& \sup_{|x|\le 1}\int_{-1}^1\left|s_n^m(x,y)\right|w(y)dy=\sup_{|x|\le 1}\int_{-1}^1\left|\sum_{r=0}^{n-1}\frac{q_{n,r}^m(x)q_{n,r}^m(y)}{\nu_{n,r}^m}\right|w(y)dy\\
&=& \sup_{|x|\le 1}\int_{-1}^1\left|\sum_{k=1}^n \Phi_{n,k}^m(x)\sum_{r=0}^{n-1}\frac{p_r(x_k^n)q_{n,r}^m(y)}{\nu_{n,r}^m}\right|w(y)dy\\
&\le& \sup_{|x|\le 1}\sum_{k=1}^n |\Phi_{n,k}^m(x)|\int_{-1}^1 \left|\sum_{r=0}^{n-1}\frac{p_r(x_k^n)q_{n,r}^m(y)}{\nu_{n,r}^m}\right|w(y)dy\\
&\le&  \C\sup_{1\le k\le n}\int_{-1}^1 \left|\sum_{r=0}^{n-1}\frac{p_r(x_k^n)q_{n,r}^m(y)}{\nu_{n,r}^m}\right|w(y)dy.
\end{eqnarray*}
Hence, to get the statement we have to prove
\begin{equation}\label{tesi-LC}
\sup_{1\le k\le n}\int_{-1}^1 \left|\sum_{r=0}^{n-1}\frac{p_r(x_k^n)q_{n,r}^m(y)}{\nu_{n,r}^m}\right|w(y)dy\le \C    
\end{equation}
where $\C>0$ denotes a constant independent of $n$.

To this aim, we observe that by \eqref{q-basis} and \eqref{p-xh}, we have
\begin{eqnarray*}
&&\sum_{r=0}^{n-1}\frac{p_r(x_k^n)q_{n,r}^m(y)}{\nu_{n,r}^m} \\
&=&
\sum_{r=0}^{n-m}p_r(x_k^n)p_r(y)+\sum_{r=n-m+1}^{n-1}\frac{p_r(x_k^n)[\mu_{n,r}^m p_r(y)-\mu_{n,2n-r}^m p_{2n-r}(y)]}{\nu_{n,r}^m}\\
&=& 
\sum_{r=0}^{n-m}p_r(x_k^n)p_r(y)+\sum_{r=n-m+1}^{n-1}\frac{\mu_{n,r}^m}{\nu_{n,r}^m}p_r(x_k^n)p_r(y)+\sum_{r=n-m+1}^{n-1}\frac{\mu_{n,2n-r}^m}{\nu_{n,r}^m}p_{2n-r}(x_h^n)p_{2n-r}(y)\\
&=& 
\sum_{r=0}^{n-m}p_r(x_k^n)p_r(y)+\sum_{r=n-m+1}^{n-1}\frac{\mu_{n,r}^m}{\nu_{n,r}^m}p_r(x_k^n)p_r(y)+\sum_{s=n+1}^{n+m-1}\frac{\mu_{n,s}^m}{\nu_{n,2n-s}^m}p_{s}(x_h^n)p_{s}(y),
\end{eqnarray*}
i.e., recalling \eqref{muj} and \eqref{nur}, we get
\[
\sum_{r=0}^{n-1}\frac{p_r(x_k^n)q_{n,r}^m(y)}{\nu_{n,r}^m}= \sum_{r=0}^{n+m-1}h_rp_r(x_k^n)p_r(y),
\]
where we set
\[
h_r:=\left\{\begin{array}{ll}
  1   &  0\le r\le n-m\\
\displaystyle   m\left(\frac{m+n-r}{m^2+(n-r)^2}\right)   &  n-m<r<n+m\\
  0 & \mbox{otherwise}.
\end{array}\right.
\]
Therefore, proving \eqref{tesi-LC} is equivalent to proving 
\[
\sup_{1\le k\le n}\int_{-1}^1 \left|\sum_{r=0}^{n+m-1}h_rp_r(x_k^n)p_r(y)\right|w(y)dy\le\C.
\]
But this follows from \cite[Lemma 5.3]{TB-GVP} whose assumptions hold since the coefficients $h_r$ satisfy
\[
\sum_{r=0}^{n+m-1} |\Delta^2 h_r|=\bigO\left( \frac 1 n \right)
\]
where $\Delta^2h_r=\Delta h_{r+1}-\Delta h_r$ and $\Delta h_r=h_{r+1}-h_r$. \Proofend

\subsection{Proof of Corollary \ref{cor-Snm}}
The uniform boundedness result
\begin{equation}\label{eq1}
\sup_n \|S_n^m\|_{L^p_w\to L^p_w} <\infty  
\end{equation}
is a direct consequence of \eqref{Sn-LC}, \eqref{Sn-norm-1} and \eqref{Sn-norm-p}.

The left--hand side inequality in \eqref{Snm-nearbest} trivially follows from $S_n^mf\in \V_n^m\subseteq \PP_{n+m-1}$ and the definition \eqref{E-best}.

To prove the right--hand side in  \eqref{Snm-nearbest}, we note that 
\begin{equation}\label{inva-S}
S_n^m P=P, \qquad \forall P\in\PP_{n-m} 
\end{equation}
follows from  $S_n^mf=f, \ \forall f\in \V_n^m\supseteq\PP_{n-m}$. 

Thus, taking $P^*\in\PP_{n-m}$ such that $\|f-P^*\|_p=E_{n-m}(f)_p$, we get 
\[
 \|S_n^mf-f\|_p \le  \|S_n^m(f-P^*)\|_p  +\|P^*-f\|_p \le (\|S_n^m\|_{L^p_w\to L^p_w} + 1) E_{n-m}(f)_p
\] 
which, due to \eqref{eq1}, concludes the proof. \Proofend

\subsection{Proof of Theorem \ref{th-LF}.} To prove \eqref{eq-LF}, we recall
that in \cite[Thm. 3.1 (f)]{TB-wave} the following inequalities are proved for the VP interpolation polynomial of any function $f$ (cf. \eqref{VP})
\begin{equation}\label{norm1-equi}
\C_1\left(\frac\pi n \sum_{i=1}^n |f(x_i^n)|\right)\le \int_{-1}^1 |V_n^mf(y)|w(y)dy
\le \C_2 \left(\frac\pi n \sum_{i=1}^n |f(x_i^n)|\right)    
\end{equation}
where $\C_1, \C_2>0$  are independent of the specific function $f$ and also independent of $n$ and $m$ in case $n/m=\bigO(1)$, which is certainly verified by our assumption $m=\lfloor \theta n\rfloor$ . 

In the case of functions $f\in \V_n^m$, since we have $V_n^m f=f$, \eqref{norm1-equi} reduces to
\[
\C_1\left(\frac\pi n \sum_{i=1}^n |f(x_i^n)|\right)\le \int_{-1}^1 |f(y)|w(y)dy
\le \C_2 \left(\frac\pi n \sum_{i=1}^n |f(x_i^n)|\right) , \quad \forall f\in \V_n^m.  
\]
Thus, if we apply \eqref{norm1-equi} to the function $f_x(y)=s_n^m(y,x)$ with  $x\in [-1,1]$ arbitrarily fixed, then the statement easily follows by taking into account that $f_x\in \V_n^m$ is a direct consequence of \eqref{sn-q} and 
$\V_n^m=\mbox{span} \{q_{n,r}^m \ : \ r=0,\ldots, n-1\}$.
\Proofend

\subsection{Proof of Theorem \ref{th-tildeSn}}
First of all,  we observe that by \eqref{Gauss} and \eqref{inva-S} we get
\begin{equation}\label{inva-Stilde}
  \tilde S_n^m P=S_n^m P=P, \qquad \forall P\in\PP_{n-m}  .
\end{equation}
So we can proceed as in the proof of Corollary \ref{cor-Snm} to prove \eqref{Stilde-nearbest} using \eqref{tildeSn-p} and $\PP_{n-m}\subseteq \V_n^m\subseteq\PP_{n+m-1}$.

Hence, let us prove \eqref{tildeSn-p}. To this aim, recalling  \eqref{eq-LF1}, Thm.~\ref{th-LC} and Cor.~\eqref{cor-Snm}, it is sufficient to prove that
we have
 \begin{equation}\label{eq-remStilde}
 \|\tilde S_n^mf\|_p\le \C \|f\|_p\qquad \mbox{with}\qquad     \C=\left\{\begin{array}{ll}
         \tilde\Lambda_n^m &  p=\infty\\
         \Lambda_n^m &  p=1\\
         (1+4\pi)\|S_n^m\|_{L^{p'}_w\to L^{p'}_w} &  1< p<\infty,
     \end{array}\right.
 \end{equation} 
with  $\frac 1{p'}+\frac 1p=1$.

In the case $p\in \{1,\infty\}$, this easily follows from the definitions that imply the following
\begin{eqnarray*}
    \|\tilde S_n^m f\|_\infty&\le& \tilde\Lambda_n^m\ \|f\|_{\ell^\infty(X_n)}\\
  \|\tilde S_n^m f\|_1&\le& \left[\sup_{1\le i\le\infty} \lambda_n^m(x_i^n)\right] \|f\|_{\ell^1(X_n)} \le \Lambda_n^m\ \|f\|_{\ell^1(X_n)} .
\end{eqnarray*}
In the case $1< p <\infty$, we use dual arguments to get \eqref{eq-remStilde}. More precisely, recalling that the dual space of $L^p_w$ is $L^{p'}_w$, with $\frac 1p+\frac 1{p'}=1$, corresponding to any $\tilde S_n^m f$ there exists a function $g\in L^{p'}$ such that $\|g\|_{p'}=1$ and we have
\begin{eqnarray*}
    \|\tilde S_n^mf\|_p&=&
    \int_{-1}^1\tilde S_n^mf(y)\ g(y)w(y)dy\\
    &=& \sum_{i=1}^n f(x_i^n) \int_{-1}^1 s_n^m(x_i^n, y)\ g(y)w(y)dy\\
    &=& \frac \pi n \sum_{i=1}^n f(x_i^n)  S_n^mg(x_i^n)\\
    &\le& \left(\frac \pi n \sum_{i=1}^n|f(x_i^n)|^p \right)^\frac 1p
    \left(\frac \pi n \sum_{i=1}^n| S_n^mg(x_i^n)|^{p'} \right)^\frac 1{p'}\\
    &=& \|f\|_{\ell^p(X_n)}\left(\frac \pi n \sum_{i=1}^n| S_n^mg(x_i^n)|^{p'} \right)^\frac 1{p'}.
    \end{eqnarray*}
On the other hand, taking into account that \cite[Thm.3.1 (f)]{TB-wave}
\[
\left(\frac \pi n \sum_{i=1}^n| F(x_i^n)|^{p'} \right)^\frac 1{p'}\le (1+4\pi) \|F\|_{p'}, \qquad \forall F\in \V_n^m,
\]
when $F=\tilde S_n^mf$, we get
\begin{eqnarray*}
    \left(\frac \pi n \sum_{i=1}^n| S_n^mg(x_i^n)|^{p'} \right)^\frac 1{p'}&\le& (1+4\pi) \|S_n^m g\|_{p'}\le (1+4\pi) \|S_n^m \|_{L^{p'}_w\to L^{p'}_w}\|g\|_{p'}\\
    &=& (1+4\pi) \|S_n^m \|_{L^{p'}_w\to L^{p'}_w} ,
\end{eqnarray*}
which, by the previous estimate, allows us to conclude that
\[
 \|\tilde S_n^mf\|_p\le (1+4\pi) \|S_n^m \|_{L^{p'}_w\to L^{p'}_w}
 \|f\|_{\ell^p(X_n)}, \qquad 1< p<\infty.
\]
\Proofend

\subsection{Proof of Proposition \ref{prop-wav}} 
Since the polynomials in \eqref{wav-ort} form an orthogonal basis of $\W_n^m$, by Def. \ref{def-wav}, we certainly have $\tilde\psi_{n,k}^m\in\W_n^m$ for any $k=1,\ldots, 2n$. Moreover, by \eqref{prod1} we have
\[
<\tilde\psi_{n,k}^m,\ \tilde\psi_{n,h}^m>_{L^2_w}=
\sum_{r=n}^{3n-1}\sigma_{r,k}^m\ \sigma_{r,h}^m, \qquad k,h=1,\ldots,2n.
\]
Hence, to complete the proof and state \eqref{prod-wav}, it is sufficient to prove the orthogonality of the square matrix $M\in \RR^{2n\times 2n}$ defined by
\[
M=[M_{r,k}]_{r,k}, \qquad
M_{r,k}=\sigma_{r,k}^m
\qquad r=n,\ldots,3n-1, \quad k=1,\ldots,2n.
\]
To this aim, in the following, we prove that 
\begin{equation}\label{tesi}
\sum_{k=1}^{2n}M_{r,k}M_{s,k}=\delta_{r,s}, \qquad r,s=n,\ldots, 3n-1.
\end{equation} 
Observe that, by \eqref{sigma}, we can write
\[
M_{r,k}=\sqrt{\frac \pi{3n}}\sigma_r(y_k^n),
\qquad r=n,\ldots,3n-1, \quad k=1,\ldots,2n
\]
where $\sigma_r$ denotes the following polynomial of degree at most $3n-1$
\[
\sigma_{r}(x)=\left\{
\begin{array}{ll}
p_n(x) & \mbox{if $r=n$}\\ [.15in]
\displaystyle \frac{p_r(x)+p_{2n-r}(x)}{\sqrt{2}} & \mbox{if $n<r<2n$}\\[.15in]
\displaystyle \frac{p_{2n}(x)+\sqrt{2}\ p_0(x)}{\sqrt{3}} &\mbox{if $r=2n$}\\ [.15in]
\displaystyle \sqrt{\frac 32} p_r(x) &\mbox{if $2n<r<3n$}\\
\end{array}\right.\qquad x\in [-1,1].
\]
Recalling that $<p_r, p_s>_{L^2_w}=\delta_{r,s}$, we easily get
\begin{equation}\label{sigma-prod}
<\sigma_r,\ \sigma_s>_{L^2_w}= \delta_{r,s}\left\{ \begin{array}{ll}
\frac 32 & 2n<r,s<3n,\\
1 & \mbox{otherwise}
\end{array}\right.\qquad 
\qquad r,s=n,\ldots, 3n-1. 
\end{equation}
On the other hand, applying the Gaussian rule on the nodes $X_{3n}=X_n\cup Y_n$, we obtain 
\begin{eqnarray*}
<\sigma_r,\ \sigma_s>_ {L^2_w}&=&
\frac \pi{3n}\sum_{k=1}^{3n}\sigma_r(x_k^{3n})\sigma_s(x_k^{3n})\\
&=&\frac\pi{3n}\left[\sum_{k=1}^{2n}\sigma_r(y_k^n)\sigma_s(y_k^n)\right]+\frac \pi{3n}\left[\sum_{k=1}^{n}\sigma_r(x_k^n)\sigma_s(x_k^n)\right]\\
&=&\sum_{k=1}^{2n}M_{r,k}M_{s,k}+\frac \pi{3n}\left[\sum_{k=1}^{n}\sigma_r(x_k^n)\sigma_s(x_k^n)\right],
\end{eqnarray*}
i.e., the following holds 
\begin{equation}\label{MMt}
\sum_{k=1}^{2n}M_{r,k}M_{s,k} = <\sigma_r,\ \sigma_s>_ {L^2_w} -
\frac \pi{3n}\left[\sum_{k=1}^{n}\sigma_r(x_k^n)\sigma_s(x_k^n)\right],\quad r,s=n,\ldots, 3n-1.
\end{equation}
By means of \eqref{sigma-prod} and \eqref{MMt}, let us prove \eqref{tesi} by distinguishing the following cases:
\begin{itemize}
\item Case $n\le s< 3n$ and $n\le r\le 2n$. Recalling \eqref{p-xh}
we easily deduce that
\[
\sigma_r(x_k^n)= 0 \qquad k=1,\ldots,n,\qquad n\le r\le 2n.
\]
Consequently,  we get
\[
\sum_{k=1}^{n}\sigma_r(x_k^n)\sigma_s(x_k^n)=0, \qquad n\le s <3n,\qquad n\le r\le 2n,
\]
and \eqref{sigma-prod}-\eqref{MMt} imply
\begin{equation}\label{sigma-1}
\sum_{k=1}^{2n}M_{r,k}M_{s,k}=<\sigma_r,\ \sigma_s>_ {L^2_w}=\delta_{r,s}, \qquad n\le s <3n,\qquad n\le r\le 2n.
\end{equation}
\item Case $2n<r<3n$ and $n\le s\le 2n$. In this case, we have
\[
\sigma_s(x_k^n)= 0 \qquad k=1,\ldots,n,\qquad n\le s\le 2n.
\]
and analogously to the previous case  we get \eqref{tesi}.
\item Case $2n<r,s<3n$. Let us set 
\[
r=3n-\ell_1, \qquad s=3n-\ell_2,\qquad\qquad \ell_1,\ell_2=1,\ldots,n-1.
\]
Recalling that
\[
\cos[(3n-\ell)t]=\cos[(2nt]\cos[(n-\ell)t]-\sin[(2nt]\sin[(n-\ell)t]
\]
we easily deduce
\[
p_{3n-\ell}(x_k^n)=-p_{n-\ell}(x_k^n), \qquad \ell=1,\ldots,n-1,\qquad k=1,\ldots,n,
\]
 which implies the following
\begin{eqnarray*}
\frac\pi n\sum_{k=1}^{n}p_r(x_k^n)p_s(x_k^n)&=&
\frac\pi n\sum_{k=1}^{n}p_{3n-\ell_1}(x_k^n)p_{3n-\ell_1}(x_k^n)
\\
&=&\frac\pi n\sum_{k=1}^{n}p_{n-\ell_1}(x_k^n)p_{n-\ell_2}(x_k^n)\\
&=&
<p_{n-\ell_1}, \ p_{n-\ell_2}>_{L^2_w}=\delta_{\ell_1, \ell_2}=\delta_{r,s},
\end{eqnarray*}
having used \eqref{Gauss} at the last line.

Consequently, for $2n<r,s<3n$ we get
\[
\frac\pi{3n} \sum_{k=1}^{n}\sigma_r(x_k^n)\sigma_s(x_k^n)=
\frac\pi{2n}  \sum_{k=1}^{n}p_r(x_k^n)p_s(x_k^n)
=\frac 12 \delta_{r,s},
\]
and by \eqref{sigma-prod}-\eqref{MMt} we conclude that
\begin{equation}\label{sigma-2}
\sum_{k=1}^{2n}M_{r,k}M_{s,k}=<\sigma_r,\ \sigma_s>_ {L^2_w} -\frac 12 \delta_{r,s}=\delta_{r,s}, \qquad 2n<r, s <3n.
\end{equation}
\end{itemize}  
\Proofend

\subsection{Proof of Theorem \ref{th-dec} }
First of all, let us observe that, for all $x\in [-1,1]$,  the following two-scale relations
\begin{eqnarray}
\label{2scale1}
    \tilde\varphi_{n,k}^m(x) &=&\sum_{j=1}^{3n}<\tilde\varphi_{n,k}^m, \ \tilde\varphi_{3n,j}^m>_ {L^2_w} \tilde\varphi_{3n,j}^m(x), \qquad k=1,\ldots,n,\\
    \label{2scale2}
  \tilde\psi_{n,h}^m(x) &=&\sum_{j=1}^{3n}<\tilde\psi_{n,h}^m, \ \tilde\varphi_{3n,j}^m>_ {L^2_w} \tilde\varphi_{3n,j}^m(x),\qquad h=1,\ldots,2n ,  
\end{eqnarray}
trivially follow by expanding $\tilde\varphi_{n,k}^m\in \V_{3n}^m$ and $\tilde\psi_{n,h}^m\in\V_{3n}^m$ in the orthonormal VP scaling basis $\{\tilde\varphi_{3n,j}^m,\ j=1,\ldots,3n\}$ of $\V_{3n}^m$.

The previous relations are governed by the rectangular matrices
\[
\begin{array}{lll}
A=[A_{k,j}] & \quad A_{k,j}=<\tilde\varphi_{n,k}^m, \ \tilde\varphi_{3n,j}^m>_ {L^2_w} & \quad k=1,\ldots, n,\qquad j=1,\ldots,3n\\
B=[B_{h,j}] &\quad  B_{h,j}=<\tilde\psi_{n,h}^m, \ \tilde\varphi_{3n,j}^m>_ {L^2_w} & \quad h=1,\ldots, 2n,\ \quad j=1,\ldots,3n.
\end{array}
\]
They form the square matrix $Q=[Q_{k,j}]\in\RR^{3n\times 3n}$ defined as follows
\[
Q=\left[\begin{array}{c}
     A \\
     B 
\end{array}\right], \qquad 
Q_{k,j}=\left\{\begin{array}{ll}
A_{k,j} & k=1,\ldots, n,\\
B_{k-n,j} & k=n+1,\ldots, 3n,
\end{array}\right. \quad j=1,\ldots, 3n.
\]
Let us prove that $Q$ is an orthogonal matrix stating that $QQ^T=I$, where $I$ denotes the identity matrix.

Indeed, taking into account the orthogonality properties stated in \eqref{compl}, Prop. \ref{prop-sca}, Prop. \ref{prop-wav},  and using \eqref{2scale1}--\eqref{2scale2}, we get
\begin{eqnarray*}
\delta_{h,k}&=&<\tilde\varphi_{n,h}^m, \ \tilde\varphi_{n,k}^m>_ {L^2_w} = \sum_{j=1}^{3n}A_{hj}A_{kj}=(AA^T)_{hk}, \qquad h,k=1,\ldots,n  \\
\delta_{h,k}&=&<\tilde\psi_{n,h}^m, \ \tilde\psi_{n,k}^m>_ {L^2_w} = \sum_{j=1}^{3n}B_{hj}B_{kj}=(BB^T)_{h,k}, \qquad h,k=1,\ldots,2n,\\  
0 &=&<\tilde\psi_{n,h}^m, \ \tilde\varphi_{n,k}^m>_ {L^2_w}=\sum_{j=1}^{3n} B_{hj} A_{kj} =(BA^T)_{hk}=(AB^T)_{kh},\qquad \begin{array}{l}
h=1,\ldots,2n,\\
k=1,\ldots,n
\end{array}\ 
\end{eqnarray*}
which implies
\begin{equation}\label{Q-ort}
QQ^T=\left[\begin{array}{c}
     A \\
     B 
\end{array}\right] \left[A^T\ B^T\right]=
\left[\begin{array}{cc}
     AA^T & AB^T \\
     BA^T & BB^T 
\end{array}\right]=I.
\end{equation}

As a consequence of \eqref{Q-ort}, we have that $Q^{-1}=Q^T$ and the two scale relations \eqref{2scale1}--\eqref{2scale2} can be inverted simply taking the transpose of the matrices $A$ and $B$. More precisely, in vector notation, setting
\[
\vec{\varphi}_n=\left[\begin{array}{c}
\tilde\varphi_{n,1}^m(x)\\
\vdots\\
\tilde\varphi_{n,n}^m(x)
\end{array}\right]\qquad\qquad
\vec{\psi}_n=\left[\begin{array}{c}
\tilde\psi_{n,1}^m(x)\\
\vdots\\
\tilde\psi_{n,2n}^m(x)
\end{array}\right]
\]
we have the following vector form of the two scale relations
\begin{equation}\label{2scale-vec}
\left[\begin{array}{c}
     \vec{\varphi}_n \\
     \vec{\psi}_n 
\end{array}\right]  =  
\left[\begin{array}{c}
     A \\
     B 
\end{array}\right]\ \vec{\varphi}_{3n} \qquad \qquad 
\vec{\varphi}_{3n} = \left[A^T\ B^T\right]\left[\begin{array}{c}
     \vec{\varphi}_n \\
     \vec{\psi}_n 
\end{array}\right] 
\end{equation}
that easily implies
\[
\left[\begin{array}{c}
     \vec{a}_n \\
     \vec{b}_n 
\end{array}\right]  =  
\left[\begin{array}{c}
     A \\
     B 
\end{array}\right]\ \vec{a}_{3n} \qquad \qquad 
\vec{a}_{3n} = \left[A^T\ B^T\right]\left[\begin{array}{c}
     \vec{a}_n \\
     \vec{b}_n 
\end{array}\right] 
\]
i.e. \eqref{dec-a} and \eqref{rec-a} hold.

Hence, to complete the proof, it remains to prove that the identities \eqref{A} and \eqref{B} hold for the entries $A_{kj}$ and $B_{hj}$.

Regarding $A_{kj}$, for $k=1,\ldots,n$ and $j=1,\ldots,3n$, 
taking into account that
\[
<q_{n,r}^m,\ q_{3n,s}^m>_ {L^2_w}=\left\{\begin{array}{lll}
  \delta_{r,s}   &  0\le r\le n-m, & 0\le s\le 3n-m \\
   \mu_{n,r}^m\delta_{r,s}-\mu_{n,2n-r}^m \delta_{2n-r,s}   &  n-m < r<n & 0\le s\le 3n-m \\
    0 &  0\le r<n & 3n-m< s<3n
\end{array}\right.
\]
 by \eqref{sca-ort}, we get
\begin{eqnarray*}
A_{kj}&=&<\tilde\varphi_{n,k}^m, \ \tilde\varphi_{3n,j}^m>_ {L^2_w}\ =\ \sum_{r=0}^{n-1}\sum_{s=0}^{3n-1}\tau_{r,k}^n\tau_{s,j}^{3n} <q_{n,r}^m, \ q_{3n,s}^m>_ {L^2_w}
\\
&=&
\sum_{r=0}^{n-m} \tau_{r,k}^n \tau_{r,j}^{3n}+
\sum_{r=n-m+1}^{n-1}\tau_{r,k}^n\left[\mu_{n,r}^m \tau_{r,j}^{3n}
-\mu_{n,2n-r}^m \tau_{2n-r,j}^{3n}\right].
\end{eqnarray*}
Similarly, with regard to $B_{hj}$, with $h=1,\ldots,2n$ and $j=1,\ldots,3n$,  taking into account that
\[
<\tilde q_{n,r}^m,\ q_{3n,s}^m>_ {L^2_w}=\left\{\begin{array}{lll}
  \mu_{n,r}^m\delta_{2n-r,s}+ \mu_{n,2n-r}^m \delta_{r,s} & n\le r<n+m & 0\le s\le 3n-m\\
  \delta_{r,s}   &  n+m\le r \le 3n-m, & 0\le s\le 3n-m \\
  0  &  3n-m < r<3n & 0\le s\le 3n-m \\
    0 &  n\le r\le 3n-m & 3n-m< s<3n\\
    v_{n,r}^m\delta_{r,s} & 3n-m< r<3n & 3n-m< s<3n
\end{array}\right.
\]
from \eqref{wav-ort2} and \eqref{sca-ort}, we deduce
\begin{eqnarray*}
B_{h,j}
&=& <\tilde\psi_{n,h}^m, \ \tilde\varphi_{3n,j}^m>_ {L^2_w}\ =\ \sum_{r=n}^{3n-1}\sum_{s=0}^{3n-1}\frac{\sigma_{r,h}^n}{\sqrt{v_{n,r}^m}}\tau_{s,j}^{3n} <\tilde q_{n,r}^m, \ q_{3n,s}^m>_ {L^2_w}\\
&=&
\sum_{r=n}^{n+m-1}\frac{\sigma_{r,h}^n}{\sqrt{v_{n,r}^m}}\left[\mu_{n,r}^m \tau_{2n-r,j}^{3n}+\mu_{n,2n-r}^m \tau_{r,j}^{3n}\right]+
\sum_{r=n+m}^{3n-m}\frac{\sigma_{r,h}^n}{\sqrt{v_{n,r}^m}} \tau_{r,j}^{3n}+\\
&& +\sum_{r=3n-m+1}^{3n-1}\frac{\sigma_{r,h}^n}{\sqrt{v_{n,r}^m}} v_{n,r}^m \tau_{r,j}^{3n}.
\end{eqnarray*}
\Proofend

\section{Conclusion}\label{sec-conclusion}
Instead of taking the fundamental VP polynomials as interpolating scaling functions to generate the corresponding VP approximation space, we define orthonormal scaling functions. Similarly, in the associated detail space, we introduce orthonormal wavelets instead of interpolating wavelets. These new bases, like the interpolatory ones, are well localized around the Chebyshev nodes to which they refer; however, they are not generated by dilation and translation of a single mother function.

For this reason, they do not require any special treatment at the boundaries, since they are naturally defined in the interval $[-1,1]$, to which any other compact interval can be mapped. Other peculiar features include the dependence on a free parameter that determines the number of vanishing moments, and a nonstandard multiresolution analysis based on a scaling factor of three rather than two.

In the paper, we derive computationally efficient decomposition and reconstruction algorithms based on the fast discrete cosine transform. Moreover, we study the approximation properties of the Fourier projection and its discrete counterpart associated with the new orthonormal scaling functions. In particular, for both approximations, we prove the uniform boundedness of the corresponding Lebesgue constants and the uniform convergence to any continuous function with near-best approximation order.

A comparison with VP interpolation reveals better performance of the new approximations in the case of less regular functions. In future work, we will present an application of the new orthonormal VP wavelets to image compression, highlighting the advantages obtained with respect to interpolating VP wavelets and other classical wavelets.
\section*{Funding}
The work of the first author was partially supported by GNCS-INDAM.
\newline
The work of the second author was partially supported by the Fund for Scientific Research–Flanders (Belgium), project
G0B0123N.

\bibliographystyle{abbrv}
\bibliography{biblioW}

\end{document}